\documentclass[preprint,12pt]{elsarticle}

\usepackage{amssymb}
\usepackage{amsmath,amsthm}

\usepackage[mathlines,pagewise]{lineno} 
\usepackage{etoolbox}          

\newcommand*\linenomathpatch[1]{%
	\cspreto{#1}{\linenomath}%
	\cspreto{#1*}{\linenomath}%
	\csappto{end#1}{\endlinenomath}%
	\csappto{end#1*}{\endlinenomath}%
}
\newcommand*\linenomathpatchAMS[1]{%
	\cspreto{#1}{\linenomathAMS}%
	\cspreto{#1*}{\linenomathAMS}%
	\csappto{end#1}{\endlinenomath}%
	\csappto{end#1*}{\endlinenomath}%
}

\expandafter\ifx\linenomath\linenomathWithnumbers
\let\linenomathAMS\linenomathWithnumbers
\patchcmd\linenomathAMS{\advance\postdisplaypenalty\linenopenalty}{}{}{}
\else
\let\linenomathAMS\linenomathNonumbers
\fi
\AtBeginDocument{%
	\linenomathpatch{equation}
	\linenomathpatchAMS{gather}
	\linenomathpatchAMS{multline}
	\linenomathpatchAMS{align}
	\linenomathpatchAMS{alignat}
	\linenomathpatchAMS{flalign}
}

\makeatletter
\patchcmd{\mmeasure@}{\measuring@true}{
	\measuring@true
	\ifnum-\linenopenaltypar>\interdisplaylinepenalty
	\advance\interdisplaylinepenalty-\linenopenalty
	\fi
}{}{}
\makeatother

\usepackage{geometry}
\usepackage{url}
\usepackage{enumitem}
\usepackage[hidelinks, colorlinks]{hyperref}
\usepackage{bookmark}
\bookmarksetup{
	numbered,
	addtohook={%
		\ifnum\bookmarkget{level}>1 %
		\bookmarksetup{numbered=false}%
		\fi
	},
}

\usepackage{framed,multirow}
\usepackage{subcaption}
\usepackage{caption}
\usepackage{graphicx}
\usepackage{multirow, multicol, makecell}
\usepackage{longtable}
\usepackage{threeparttable}

\usepackage{enumitem}

\usepackage{stackengine}
\newcommand{\ubar}[1]{\stackunder[1.2pt]{$#1$}{\rule{1.0ex}{.075ex}}}

\usepackage{tikz}
\usetikzlibrary {arrows.meta}
\usetikzlibrary{shapes.geometric, arrows}
\tikzstyle{startstop} = [rectangle, rounded corners, text width=3cm,
minimum width=3.05cm, 
minimum height=0.5cm,
text centered, 
draw=black, 
fill=portlandorange!10]

\tikzstyle{io} = [trapezium, rounded corners, text width=3cm,
trapezium stretches=true, 
trapezium left angle=70, 
trapezium right angle=110, 
minimum width=3.05cm, 
minimum height=0.5cm, text centered, 
draw=black, fill=screaminggreen!30]

\tikzstyle{process} = [rectangle, rounded corners, text width=3cm,
minimum width=3.05cm, 
minimum height=0.5cm, 
text centered, 
text width=3.05cm, 
draw=black, 
fill=pistachio!10]

\tikzstyle{decision} = [diamond, aspect=1.3,
minimum width=3.05cm, 
minimum height=1cm, 
text centered, 
draw=black, 
fill=pistachio!10]
\tikzstyle{arrow} = [thin,->,-{Stealth}]

\definecolor{portlandorange}{rgb}{1.0, 0.35, 0.21}
\definecolor{pistachio}{rgb}{0.58, 0.77, 0.45}
\definecolor{screaminggreen}{rgb}{0.46, 1.0, 0.44}
\definecolor{gold}{rgb}{0.85, 0.65, 0.13}

\usepackage{cleveref}

\usepackage{framed,multirow}

\usepackage{amssymb}

\usepackage{amsfonts}
\usepackage{mathrsfs}
\usepackage{upgreek}

\usepackage{graphicx}
\usepackage{epstopdf}
\usepackage{algorithmic}
\ifpdf
\DeclareGraphicsExtensions{.eps,.pdf,.png,.jpg}
\else
\DeclareGraphicsExtensions{.eps}
\fi
\usepackage[export]{adjustbox}

\usepackage{xspace}
\usepackage{bm}
\usepackage{physics}
\usepackage{stmaryrd}
\usepackage{cleveref}
\usepackage[labelfont=bf,textfont=normalfont]{caption}
\usepackage{subcaption}
\usepackage{mathtools}

\usepackage{booktabs}
\usepackage{multirow, multicol, makecell}
\usepackage{longtable}
\usepackage{threeparttable}

\crefname{assumption}{Assumption}{Assumptions}
\crefname{definition}{Definition}{Definitions}
\crefname{lemma}{Lemma}{Lemmas}
\crefname{remark}{Remark}{Remarks}
\crefname{theorem}{Theorem}{Theorems}
\crefname{proposition}{Proposition}{Propositions}
\crefname{section}{Section}{Sections}
\crefname{figure}{Fig.}{Figs.}
\crefname{equation}{}{}
\crefname{table}{Table}{Tables}
\crefname{appendix}{}{}

\newcommand{\E}[0]{\mathsf{E}}

\newcommand{\T}[0]{\mathsf{T}}

\newcommand{\A}[0]{\mathsf{A}}

\newcommand{\J}[0]{\mathsf{J}}

\newcommand{\V}[0]{\mathsf{V}}

\newcommand{\etal}[0]{{et~al.\@}\xspace}
\newcommand{\eg}[0]{{e.g.\@}\xspace}
\newcommand{\ie}[0]{{i.e.\@}\xspace}

\newcommand{\ignore}[1]{} 

\newcommand{\IR}[1]{\mathbb{R}^{#1}}%
\newcommand{\IRtwo}[2]{\mathbb{R}^{{#1}\times{#2}}}%

\newcommand{\pder[2]}[]{\dfrac{\partial#1}{\partial#2}}%
\global\long\def\norm#1{\left\vert \left\vert #1\right\vert \right\vert }%
\global\long\def\fn#1{\mathcal{#1}}%
\global\long\def\mathds#1{\mathds{#1}}%

\makeatletter
\let\overlinewithoriginalheight\overline
\newcommand*\overlinewithlessheight[1]{{\mathpalette\overline@aux{#1}}}
\newcommand*\overline@aux[2]{
	\begingroup
	\count0=\fam 
	\setbox0=\hbox{$\m@th #1\fam=\count0 #2$}
	\@tempdima=.4\ht0
	\setbox0=\hbox{$\m@th #1\fam=\count0\overlinewithoriginalheight{#2}$}%
	\advance\@tempdima by .6\ht0
	\ht0=\@tempdima 
	\usebox0
	\endgroup%
}
\let\overline\overlinewithlessheight

\journal{Computers \& Mathematics with Applications}

\makeatletter
\nopreprintlinetrue

\def\ps@pprintTitle{%
    \let\@oddhead\@empty
    \let\@evenhead\@empty
    \def\@oddfoot{\hfil\thepage\hfil} 
    \let\@evenfoot\@oddfoot%
}
\makeatother
\begin{document}

\begin{frontmatter}

\title{Very high-order symmetric positive-interior quadrature rules on triangles and tetrahedra}

\author[1]{Zelalem Arega Worku\corref{cor1}}
\cortext[cor1]{Corresponding author: } 
\ead{zelalem.worku@mail.utoronto.ca}

\author[2]{Jason E. Hicken}
\ead{hickej2@rpi.edu}

\author[1]{David W. Zingg}
\ead{dwz@oddjob.utias.utoronto.ca}

\address[1]{Institute for Aerospace Studies, University of Toronto, 
	4925 Dufferin St, Toronto, ON, 
	M3H 5T6,
	Canada}

\address[2]{Department of Mechanical, Aerospace, and Nuclear Engineering, Rensselaer Polytechnic Institute,
	110 8th St, Troy, NY,
	2180-3590, 
	USA}

\begin{abstract}
We present novel fully-symmetric quadrature rules with positive weights and strictly interior nodes of degrees up to 84 on triangles and 40 on tetrahedra. Initial guesses for solving the nonlinear systems of equations needed to derive quadrature rules are generated by forming tensor-product structures on quadrilateral/hexahedral subdomains of the simplices using the Legendre-Gauss nodes on the first half of the line reference element. In combination with a methodology for node elimination, these initial guesses lead to the development of highly efficient quadrature rules, even for very high polynomial degrees. Using existing estimates of the minimum number of quadrature points for a given degree, we show that the derived quadrature rules on triangles and tetrahedra are more than 95\% and 80\% efficient, respectively, for almost all degrees. The accuracy of the quadrature rules is demonstrated through numerical examples.
\end{abstract}

\begin{keyword}
Positive-interior quadrature \sep Gaussian quadrature \sep Numerical integration \sep Triangle \sep Tetrahedron 

\MSC[2020] 65D32 \sep 65M06 \sep 65M60
\end{keyword}

\end{frontmatter}


\section{Introduction}
Quadrature rules are indispensible tools in numerical methods for solving partial differential equations (PDEs), providing an efficient means of approximating integrals. Often numerical discretizations of problems involving complex geometries use simplicial meshes, necessitating the development of quadrature rules on triangles and tetrahedra. While several studies have proposed efficient high-order quadrature rules on these elements, \eg, \cite{dunavant1985high,zhang2009set,xiao2010numerical,witherden2015identification,jaskowiec2021high,chuluunbaatar2022pi}, futher improvements in terms of polynomial order and efficiency remain potentially beneficial. Currently, the highest degree symmetric quadrature rules with positive weights and interior nodes (PI) on triangles and tetrahedra are 50 \cite{xiao2010numerical} and 20 \cite{jaskowiec2021high,chuluunbaatar2022pi}, respectively. These accuracy levels, however, may not be sufficiently high in scenarios where a degree $p$ finite-element or discontinuous Galerkin (DG) discretization employs very high-order quadrature rules, \eg, $3p$ up to $6p$ \cite{kirby2003aliasing,persson2009discontinuous,williams2019analysis}, for accuracy or stability purposes. In general, the development of very high-order discretization operators is limited by the availability of sufficiently accurate quadrature rules. For example, the construction of summation-by-parts (SBP) operators of degree $p$, which are important for the design of provably stable discretizations of PDEs \cite{fernandez2014review,svard2014review}, requires at least degree $2p-1$ quadrature rules \cite{hicken2016multidimensional,fernandez2018simultaneous}; consequently, these operators can currently be constructed only up to degree 10 on tetrahedra. The development of very high-order quadrature rules can significantly improve the design of numerical methods for PDEs, with their impact potentially extending to the many other fields where quadrature rules are applied.

The properties of the quadrature rule used in a numerical method can affect the accuracy, efficiency, and stability of the entire simulation. Symmetry of the quadrature rules ensures that no unphysical asymmetry is introduced to a problem of interest \cite{witherden2014analysis,williams2020family}. Furthermore, symmetric nodes lead to difference operators with inherent symmetry, which can be exploited to reduce memory requirements. Positivity of the weights guarantees that the quadrature approximations of integrals of positive functions remain positive and prevents errors due to cancellations of quadrature contributions at each node. Furthermore, positivity of the weights is essential in the design of stable numerical methods, for example in the context of SBP discretizations, and it produces a valid norm that can be used to bound the energy or entropy functions associated with PDEs. Finally, the constraint that all quadrature nodes must be strictly inside the domain can be essential for practical implementations and approximation of integrals of functions with singularities at element facets.

Gaussian-type PI quadrature rules have been studied extensively, and it is well-known that Gaussian quadrature rules are optimal in one dimension \cite{stroud1966gaussian}. In some cases, quadrature rules with a subset of nodes placed at the facets, such as those presented in \cite{worku2024quadrature}, can be more efficient for numerical discretizations of PDEs despite having more nodes than PI quadrature rules. In many cases, however, the small number of nodes of PI quadrature rules translates into higher efficiency. While very efficient symmetric PI rules, up to degree 50, have been constructed on triangles \cite{xiao2010numerical}, the highest available degree on the tetrahedron has been limited to 15 \cite{xiao2010numerical,zhang2009set} until the recent extension up to degree 20 \cite{jaskowiec2021high,chuluunbaatar2022pi}. This is because the construction of quadrature rules in three dimensions is notoriously difficult for large values of $p$ \cite{solin2003higher,jaskowiec2020high}. One of the primary challenges arises from the absence of good-quality initial guesses, which are crucial for the convergence of nonlinear equation solvers, particularly as the polynomial degree increases. In this paper, we present a novel approach for generating initial guesses to solve the nonlinear equations needed to derive quadrature rules. This approach involves forming tensor-product structures on quadrilateral/hexahedral subdomains of the triangular/tetrahedral reference elements using the Legendre-Gauss (LG) nodes on the first half of the line reference element. The initial guesses, coupled with a methodology for node elimination, are used to derive efficient fully-symmetric PI quadrature rules of degree up to $84$ on triangles and $40$ on tetrahedra. 

The remaining sections are organized as follows. \cref{sec:pre} presents the symmetry orbits on the reference triangle and tetrahedron, the systems of nonlinear equations that need to be solved, and some lower-bound estimates of Gaussian-type quadrature rules from the literature. \cref{sec:methodology} delineates the general algorithm for deriving the quadrature rules and the method of solution. The initial guess generation is presented in \cref{sec:initial guess}, which is followed by a discussion on the derivation and efficiency of the line-Legendre-Gauss (line-LG) quadrature rules in \cref{sec:Line-LG}. \cref{sec:node elim} presents a node elimination strategy. A discussion on the derived novel quadrature rules is provided in \cref{sec:new quad}. Finally, some numerical results are presented in \cref{sec:numerical}, followed by conclusions in \cref{sec:conclusions}.

\section{Preliminaries and Background} \label{sec:pre}
For the construction of the quadrature rules in this work, we consider the triangle and tetrahedron reference elements defined, respectively, as 
\begin{align} \label{eq:ref elems}
	&\Omega_{\text{tri}}=\{\left(x_{1},x_{2}\right)\mid x_{1},x_{2}\ge-1; \; x_{1}+x_{2}\le0\},\\
	&\Omega_{\text{tet}}=\{\left(x_{1},x_{2},x_{3}\right)\mid x_{1},x_{2},x_{3}\ge-1; \; x_{1}+x_{2}+x_{3}\le-1\}, 
\end{align}
where $\bm{x}=[x_{1},\dots,x_{d}]^T$ denotes the Cartesian coordinates on the reference elements. The line reference element is used in the initial guess generation procedure and is defined as $\Omega_{\text{line}}=\{x_{1} \mid -1\le x_{1}\le 1\} $.

There are three symmetry groups on the triangle, and five on the tetrahedron \cite{felippa2004compendium,zhang2009set,witherden2015identification}. On the triangle, the symmetry groups, in barycentric coordinates, are permutations of
\begin{equation}\label{eq:sym tri}
    \begin{aligned}
        S_{1} &= \left(\frac{1}{3}, \frac{1}{3}, \frac{1}{3}\right), &&
        S_{21} = (\alpha,\alpha,1-2\alpha),&&
        S_{111}  = (\alpha,\beta,1-\alpha-\beta),
    \end{aligned}
\end{equation}
where $ \alpha $ and $ \beta $ are parameters such that the quadrature points lie strictly in the interior of the domain. Similarly, the symmetry groups on the reference tetrahedron are permutations of
\begin{equation}\label{eq:sym tet}
	\medmuskip=-0mu
    \begin{aligned}
        S_{1} &= \left(\frac{1}{4},\frac{1}{4},\frac{1}{4},\frac{1}{4}\right),
        &S_{31} &= \left(\alpha,\alpha,\alpha,1-3\alpha\right),
        &S_{22} &= \left(\alpha,\alpha,\frac{1}{2}-\alpha,\frac{1}{2}-\alpha\right),\\
        S_{211} &= \left(\alpha,\alpha,\beta,1-2\alpha-\beta\right),
        &S_{1111} &= \left(\alpha,\beta,\gamma,1-\alpha-\beta-\gamma\right). &&
    \end{aligned}
\end{equation}

The problem statement for deriving a quadrature rule on the domain $\Omega$ can be posed as: find $ \bm{y} $ and $ \bm{w} $ such that 
\begin{equation}\label{eq:prob}
    \int_{\Omega}\fn P_{j}\left({\bm{x}}\right)\dd{\Omega}=\sum_{i=1}^{n_q}{w}_{i}\fn P_{j}\left(\bm{y}_{i}\right),\qquad j\in\{1,\dots,n_{b}\},
\end{equation}
where $\bm{y}$ is the vector of the coordinate tuples of the nodes, $\bm{w}$ is the vector of quadrature weights, and  $ n_q $ and $ n_{b} $ denote the number of quadrature points and polynomial basis functions, respectively. For a degree $ q $ accurate quadrature rule on a simplex, there are $ n_{b}= {{q+d} \choose{d}} $ polynomial basis functions. Rewriting the problem statement in a matrix form, we have
\begin{equation}\label{eq:prob matrix}
    \bm{g} \coloneqq \V^{T}\bm{w}-\bm{f}=\bm{0},
\end{equation}
where $ \V $ is the Vandermonde matrix containing evaluations of the basis functions at each node along its columns and $ \bm{f}=\left[\int_{\Omega}\fn P_{1}\left({\bm{x}}\right)\dd{\Omega},\dots,\int_{\Omega}\fn P_{n_{b}}\left({\bm{x}}\right)\dd{\Omega}\right]^{T} $. We use the orthonormal Proriol-Koornwinder-Dubiner (PKD) \cite{proriol1957family,koornwinder1975two,dubiner1991spectral} basis functions, which produce well-conditioned Vandermonde matrices and are convenient, as all except the first entry of $ \bm{f} $ are zero due to the orthogonality of the basis functions. 

\subsection{Lower-Bound Estimates for the PI Quadrature Rules}\label{sec:estimate}
It is well known that the one-dimensional LG quadrature rules are optimal in the sense that they use the smallest possible number of nodes to integrate polynomials exactly \cite{stroud1966gaussian}. Such optimal rules in higher dimensions are not known, but there are some results that estimate the lower bounds of the number of nodes required to integrate polynomials exactly. In this work, we are interested not only in the estimates of the total number of nodes but also in the number of nodes in each symmetry orbit, which serves as a guide for node elimination. For the triangle, we use the estimates of Lyness and Jespersen \cite{lyness1975moderate}, which can be stated as follows. Let $\alpha_s = [3,-4,-1,0,-1,-4]$, $\alpha_q = \alpha_s[(q\mod 6)+1]$, and $\fn{E}(q,\alpha_q) = ((q+3)^2+\alpha_q)/12$; then the lower-bound estimate for each symmetry group is given by
\begin{equation}
	\begin{aligned}\label{eq:tri bound}
		|S_{111}|=&\begin{cases}
			0, & q<6\\
			\left\lfloor \frac{\fn E\left(q-6,\alpha_{q}\right)+2}{3}\right\rfloor , & q\ge6
		\end{cases},&|S_{12}|&=\left\lfloor \frac{\fn E\left(q,\alpha_{q}\right)-3S_{111}}{2}\right\rfloor ,\\|S_{1}|=&\begin{cases}
			0, & 1+2S_{12}+3S_{111}>\fn E\left(q,\alpha_{q}\right),\\
			1, & \text{otherwise},
		\end{cases}&&
	\end{aligned}
\end{equation}
where $ \lfloor\cdot\rfloor$ indicates the floor operator. The total number of nodes on the triangle is computed as 
\begin{equation}
	n_q= |S_{1}| + 3|S_{21}| + 6|S_{111}|.
\end{equation}

For the tetrahedron, we use the estimates of Wang and Papanicolopulos \cite{wang2023explicit}, which can be stated as follows. Let $\fn{I}(q,a)$ be the function
\begin{equation}
	\fn I\left(q,a\right)=\begin{cases}
		1, & q\ge a,\\
		0, & \text{otherwise},
	\end{cases}
\end{equation}
and define the following coefficients
\begin{equation}
	\begin{aligned}
		q_{12}&=q-12,\qquad m_{2}=\left\lfloor \frac{q}{2}-1\right\rfloor \fn I(q,4),\qquad m_{3}=\left\lfloor \left(\frac{q}{2}-2\right)^{2}\right\rfloor \fn I(q,6),\\
		m_{4}&=\bigg{\lfloor}\frac{\left(q_{12}+4\right)^{3}+3\left(q_{12}+4\right)^{2}-9\left(q_{12}+4\right)\left((q_{12}+4)\bmod2\right)}{144}\bigg{\rceil}\;\fn I(q_{12},0),\\
		m_{e}&=\bigg{\lfloor}\frac{\left(q+4\right)^{3}+3\left(q+4\right)^{2}-9\left(q+4\right)\left(\left(q+4\right)\bmod2\right)}{144}\bigg{\rceil}\;\fn I(q,0),
	\end{aligned}
\end{equation}
where $\lfloor \cdot \rceil$ denotes rounding to the closest integer. Then a quadrature rule on the tetrahedron with the following number of nodes in each symmetry group provides a lower-bound estimate of the number of nodes required to satisfy a degree $q$ quadrature accuracy, 
\begin{equation}\label{eq:tet bound}
	\medmuskip=-0mu
	\begin{aligned}
		|S_{1111}|&\!=\! \left\lceil \frac{m_{4}}{4}\right\rceil \!, \;\; |S_{211}| \!=\! \left\lceil \frac{m_{4}+m_{3}-4S_{1111}}{3}\right\rceil \!, \;\; |S_{22}|\!=\!\left\lceil \frac{m_{4}+m_{3}+m_{2}-3S_{211}-4S_{1111}}{2}\right\rceil \!,
		\\|S_{31}|&\!=\! \left\lfloor \frac{m_{e}-2S_{22}-3S_{211}-4S_{1111}}{2}\right\rfloor \!,\;\;|S_{1}| \!=\! m_{e}-2S_{31}-2S_{22}-3S_{211}-4S_{1111},
	\end{aligned}
\end{equation}
where $ \lceil \cdot \rceil$ denotes the ceiling operator. The total number of nodes on the tetrahedron is computed as
\begin{equation}
	n_q  = |S_{1}| + 4|S_{31}| + 6|S_{22}| + 12|S_{211}| + 24|S_{1111}|. 
\end{equation}

We note that the lower-bound estimates stated above may not produce optimal quadrature rules as the assumptions used to derive them are not necessarily true for all cases. Furthermore, the quadrature rules satisfying the estimates may require placing nodes outside the domain or using negative weights; we refer interested readers to \cite{lyness1975moderate} and \cite{wang2023explicit} for further details. With these caveats aside, we use these lower-bound estimates to define efficiency of the PI quadrature rules conservatively as  
\begin{equation} \label{eq:eff quad}
	e_{q} \coloneqq \frac{n_q^m}{n_q},
\end{equation}
where $n_q$ is the number of nodes of a given degree $q$ quadrature rule and $n_q^m$ is the number of nodes given by the lower-bound estimates stated in \cref{eq:tri bound} for triangles and \cref{eq:tet bound} for tetrahedra.

\section{Methodology}\label{sec:methodology}
The first step in solving the quadrature problem, \cref{eq:prob matrix}, is generating some intial guesses for the quadrature points and weights, which are used to compute the Vandermonde matrix and weight vector. In general, the initial guesses do not satisfy \cref{eq:prob matrix}, and will be used to iteratively reduce the residual error to machine precision. Generating initial guesses that fall within the radius of convergence of the (gradient-based) iterative solver is a challenging task and a key factor limiting the derivation of very high-order quadrature rules on tetrahedra. 

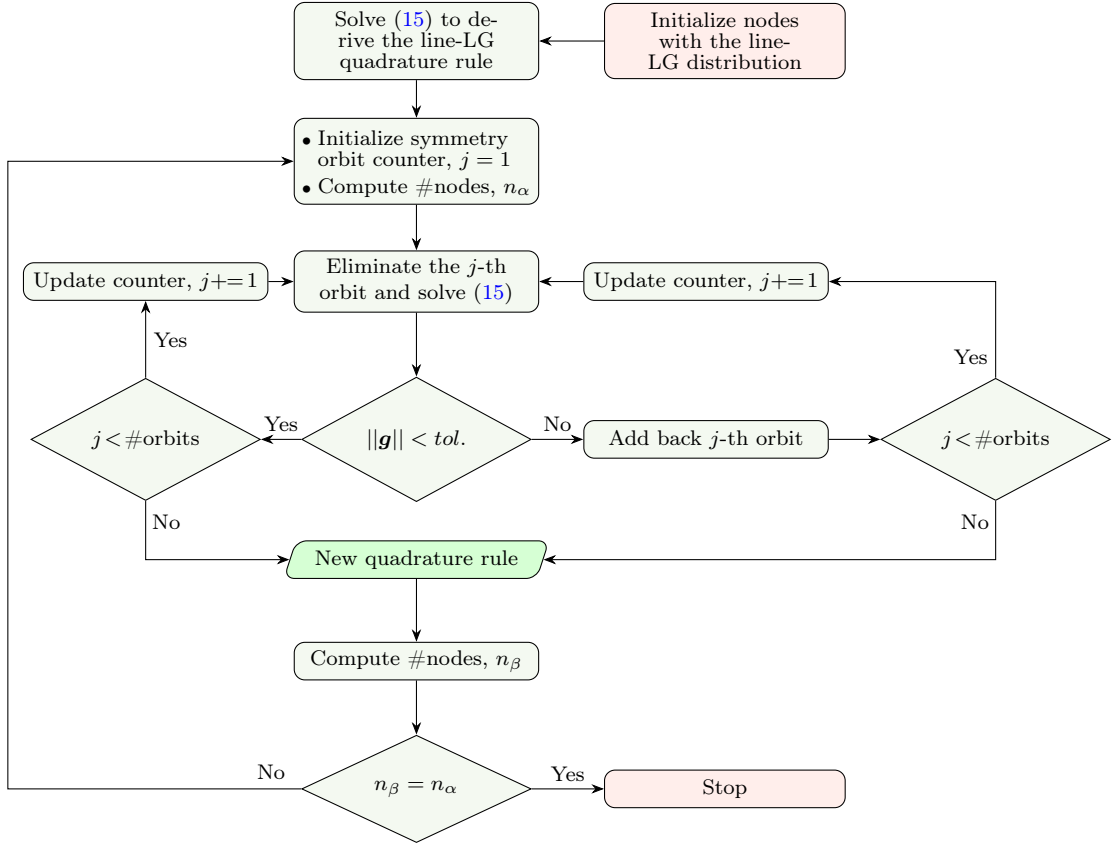
\begin{figure}[t]
	\centering
	\begin{tikzpicture}[node distance=1.6cm]
		{\fontsize{8pt}{10pt}\selectfont
		\node (LGquad) [process] {
			\baselineskip=8pt
			Solve \cref{eq:min prob} to derive the line-LG quadrature rule \par};
		\node (start) [startstop, minimum width=3cm, minimum height=1cm, text centered, text width=3cm, draw=black, right of=LGquad, xshift=2.5cm] {\baselineskip=8pt Initialize nodes with the line-LG distribution \par};
		\node (setcnt) [process, below of=LGquad, align=left] {
			\baselineskip=8pt
			\setlength\labelsep {\dimexpr\labelsep - 0.40em\relax}  
			\begin{itemize}[nosep,itemsep=2pt,leftmargin=*, label={\tiny\textbullet}]
				\vspace{-0.025cm}
				\item Initialize symmetry orbit counter, $j=1$
				\item Compute \#nodes, $n_{\alpha}$
			\end{itemize}\par};
		\node (elim) [process, below of=setcnt]{
			Eliminate the $j$-th orbit and solve \cref{eq:min prob}};
		
		\node (dec1) [decision, below of=elim, yshift=-0.5cm] {$\norm{\bm{g}}<tol.$};
		\node (dec2a) [decision, left of=dec1, xshift=-2cm] {$\! j \!<\!\text{\#orbits}\!$};
		\node (update_left) [process, left of=elim, xshift=-2cm]{Update counter, $j +\!\!=\! 1$};
		\node (pro2b) [process, right of=dec1, xshift=2.25cm] {Add back $j$-th orbit};
		\node (dec2b) [decision, right of=pro2b, xshift=2.25cm] {$\! j \!<\!\text{\#orbits}\!$};
		\node (update_right) [process, right of=elim, xshift=2.25cm]{Update counter, $j +\!\!=\! 1$};
		
		\node (out1) [io, below of=dec1] {New quadrature rule};
		\node (nnode)[process, below of=out1, yshift=0.25cm]{Compute \#nodes, $n_{\beta}$};
		\node (dec3) [decision, below of=nnode, yshift=-0.1cm] {$n_{\beta}=n_{\alpha}$};
		\node (stop) [startstop, right of=dec3, xshift=2.5cm] {Stop};
		
		\draw [arrow] (start) -- (LGquad);
		\draw [arrow] (LGquad) -- (setcnt);
		\draw [arrow] (setcnt) -- (elim);
		\draw [arrow] (elim) -- (dec1);
		\draw [arrow] (dec1) -- node[anchor=south] {Yes} (dec2a);
		\draw [arrow] (dec1) -- node[anchor=south] {No} (pro2b);
		\draw [arrow] (pro2b) -- (dec2b);
		\draw [arrow] (dec2a) |- node[xshift=0.275cm, yshift=0.5cm]{No}(out1);
		\draw [arrow] (dec2b) |- node[xshift=-0.275cm, yshift=0.5cm]{No}(out1);
		\draw [arrow] (out1) -- (nnode);
		\draw [arrow] (dec2a) -- node[anchor=west] {Yes} (update_left);
		\draw [arrow] (update_left) -- (elim);
		\draw [arrow] (dec2b) |- node[anchor=east, yshift=-1cm] {Yes} (update_right);
		\draw [arrow] (update_right) -- (elim);
		\draw [arrow] (nnode) -- (dec3);
		\draw [arrow] (dec3.west) -- ++(-3.9,0)  |- (setcnt) node[at start, xshift=3.5cm, yshift=0.225cm] {No};
		\draw [arrow] (dec3) -- node[anchor=south]{Yes}(stop);
		}
	\end{tikzpicture}
	\caption{\label{fig:flowchart}Algorithm flowchart for generating symmetric PI quadrature rules.}
\end{figure}

The methodology employed in this study can be summarized as follows. First, we introduce an efficient approach to generate initial guesses that is based on mapping the LG quadrature nodes in the first half of the line reference element onto specific lines in the simplex. Next, we set \cref{eq:prob matrix} as a minimization problem and solve it using a gradient-based iterative solver to obtain intermediate quadrature rules which we refer to as line-LG quadrature rules. Finally, we apply a node elimination technique to the line-LG rules and derive novel quadrature rules with improved efficiency. The quadrature rules obtained after node elimination will be referred to as ``new'' in tables and figures. Further details of the algorithm are illustrated in \cref{fig:flowchart}. 

The method of solution used in this work is similar to that described in \cite{worku2024quadrature} for quadrature rules with a subset of nodes placed at the facets for the construction of SBP diagonal-$\E$ multidimensional SBP operators \cite{chen2017entropy,hicken2016multidimensional}. In this work, the approach is applied to derive quadrature rules with strictly interior nodes, using only a gradient-based iterative solver instead of the coupled gradient and stochastic-based solver employed in \cite{worku2024quadrature}. For completeness, we describe the approach with some adjustments to accommodate the aforementioned differences. 

Using the initial guesses, it is possible to compute the coordinates of the $ i$-th node by applying the transformation,
\begin{equation} \label{eq:x=T lambda}
    \bm{y}_{i} = \T^{T} \bm{\lambda}_{k}^{(i)},
\end{equation}
where $ \T \in \IRtwo{(d+1)}{d}$ contains the coordinates of the $ d+1 $ vertices in its rows and $ \bm{\lambda}_{k}^{(i)} $ is the $ k $-th permutation of the barycentric coordinates of the symmetry group that corresponds to the $ i $-th node. The weight vector, $ \bm{w} $, is constructed by assigning equal weights to all nodes in the same symmetry group. 

We note that the problem in \cref{eq:prob matrix} can be written equivalently as a minimization problem,
\begin{equation}\label{eq:min prob}
    \min_{\bm{\tau}}{\frac{1}{2}\bm{g}^T\bm{g}},
\end{equation}
where $ \bm{\tau} = [\widehat{\bm{\lambda}}^T,\widehat{\bm{w}}^T]^T $, and $\widehat{\bm{\lambda}} $ and $\widehat{\bm{w}}$  are vectors of all the parameters and weights associated with each symmetry group.  The  Levenberg-Marquardt (LM) algorithm \cite{levenberg1944method,marquardt1963algorithm} is used to solve \cref{eq:min prob}. The LM algorithm computes the step direction, $ \bm{h} $ as  
\begin{equation}\label{eq:step}
 {\bm{h}} =-{\A}^{+} {\J^T \bm{g}},	
\end{equation}
where $ \A = \J^T \J + \nu \text{diag}(\J^T\J )$, $ \nu > 0$ is a parameter that controls the scale of exploration, and $ (\cdot)^{+} $ denotes the Moore-Penrose pseudo-inverse. The matrix $ \J\in \IRtwo{n_{b}}{n_{\tau}} $ is the Jacobian matrix given by
\begin{align}
    \J_{(i,j)} &= \pder[\bm{g}_{i}]{\bm{\tau}_{j}},
\end{align}
and $ n_{\tau} $ is the sum of the number of parameters and weights. The Jacobian can also be written in terms of block matrices as
\begin{align}
    \J=\left[\sum_{k=1}^{d}\V_{x_{k}}^{T}\text{diag}(\bm{w})\pder[\bm{y}_{(:,k)}]{\widehat{\bm{\lambda}}},\V^{T}\pder[\bm{w}]{\widehat{\bm{w}}}\right],
\end{align}
where $ \V_{{x}_k} $ is the $ k $-th direction derivative of $ \V $ and $  \bm{y}_{(:,k)}  $is the $ k $-th direction component vector of $ \bm{y} $. The matrix $ \partial{\bm{y}_{(:,k)}}/\partial{\widehat{\bm{\lambda}}} $ is computed using the relation in \cref{eq:x=T lambda}, and $ \partial{\bm{w}}/\partial{\widehat{\bm{w}}} $ is a matrix of zeros and ones.

The algorithm starts with an initial guess, $ \bm{\tau}^{(0)} $, and the value of $ \bm{\tau} $ at the $ n$-th iteration is updated as 
\begin{equation}\label{eq:tau update}
    \bm{\tau}^{(n+1)} = \bm{\tau}^{(n)} + \eta^{(n)}\bm{h}^{(n)},
\end{equation}
where $ \eta^{(n)} = 1$ is used unless a negative weight is encountered. If a negative weight is encountered at the $ i $-th entry of $ \bm{\tau}^{(n+1)} $, then the update is recomputed using
\begin{equation}
    \eta^{(n)} = \frac{(\varepsilon - \bm{\tau}^{(n)}_{i})}{\bm{h}^{(n)}_{i}},
\end{equation}
where $ \varepsilon > 0$ is an arbitrary lower bound for the update of the negative weight, \eg, we use $ \varepsilon = 10^{-4} $. All of the quadrature rules in this work are obtained using the open-source Julia code \texttt{SummationByParts.jl }\cite{hicken2023summationbyparts}.

\section{Initial Guess Generation}\label{sec:initial guess}
\begin{figure}[t!]
	\centering
	\begin{subfigure}{0.30\textwidth}
		\centering
		\includegraphics[scale=0.26]{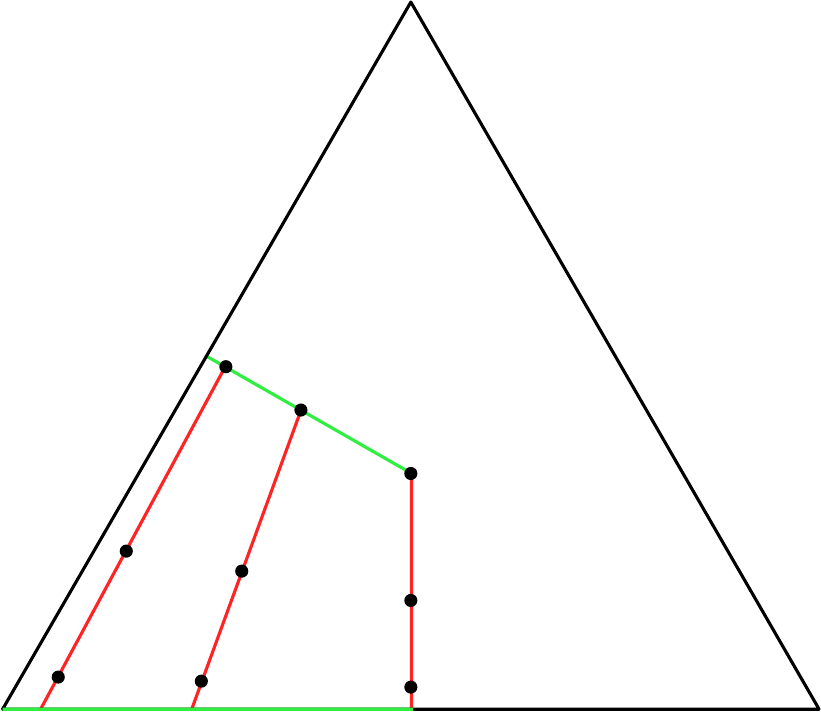}
		\caption{\label{fig:lg line tri} Half line LG  node mappings}
	\end{subfigure}\hfill
	\begin{subfigure}{0.30\textwidth}
		\centering
		\includegraphics[scale=0.26]{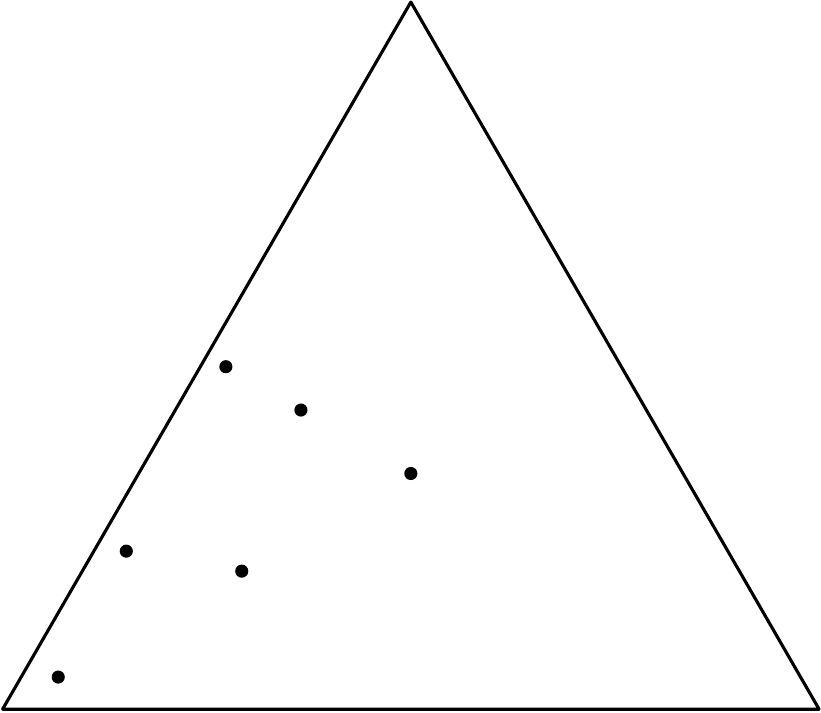}
		\caption{\label{fig:unique sym tri} Unique symmetry orbits}
	\end{subfigure}\hfill
	\begin{subfigure}{0.30\textwidth}
		\centering
		\includegraphics[scale=0.26]{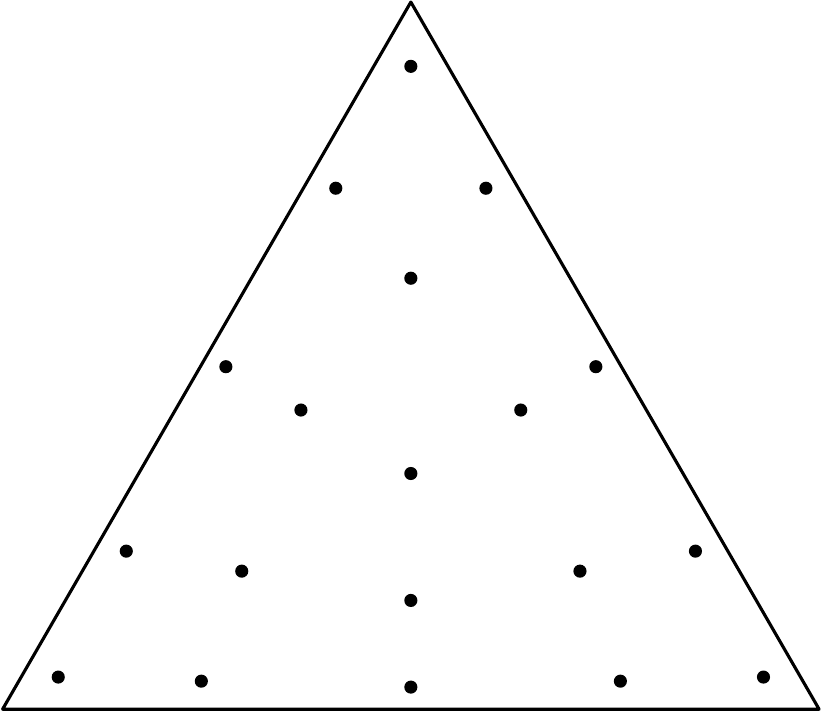}
		\caption{\label{fig:init tri} Initial node distribution}
	\end{subfigure}
	\\
	\begin{subfigure}{0.30\textwidth}
		\centering
		\includegraphics[scale=0.0130]{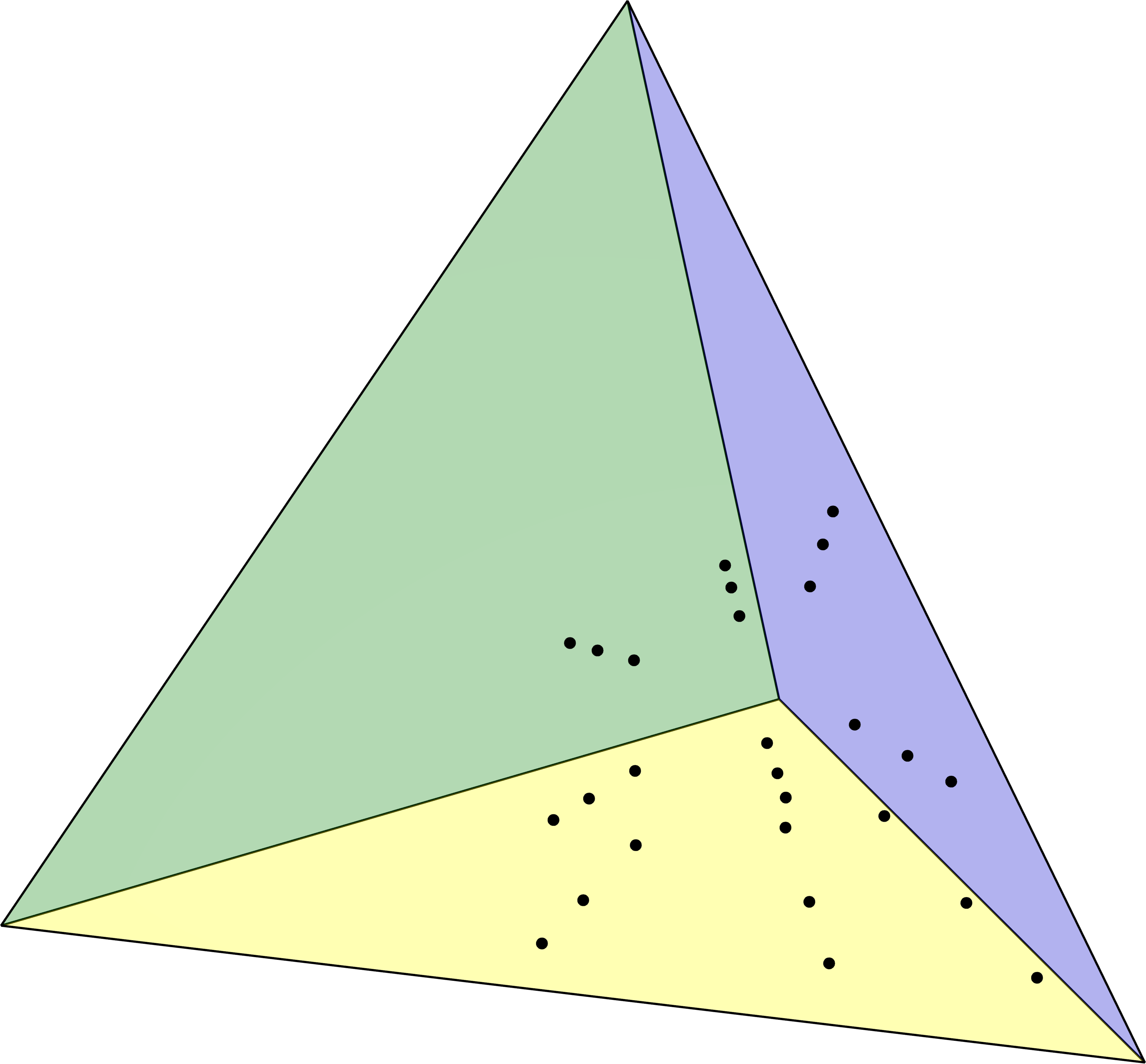}
		\caption{\label{fig:lg line tet} Half line LG  node mappings}
	\end{subfigure}\hfill
	\begin{subfigure}{0.30\textwidth}
		\centering
		\includegraphics[scale=0.0130]{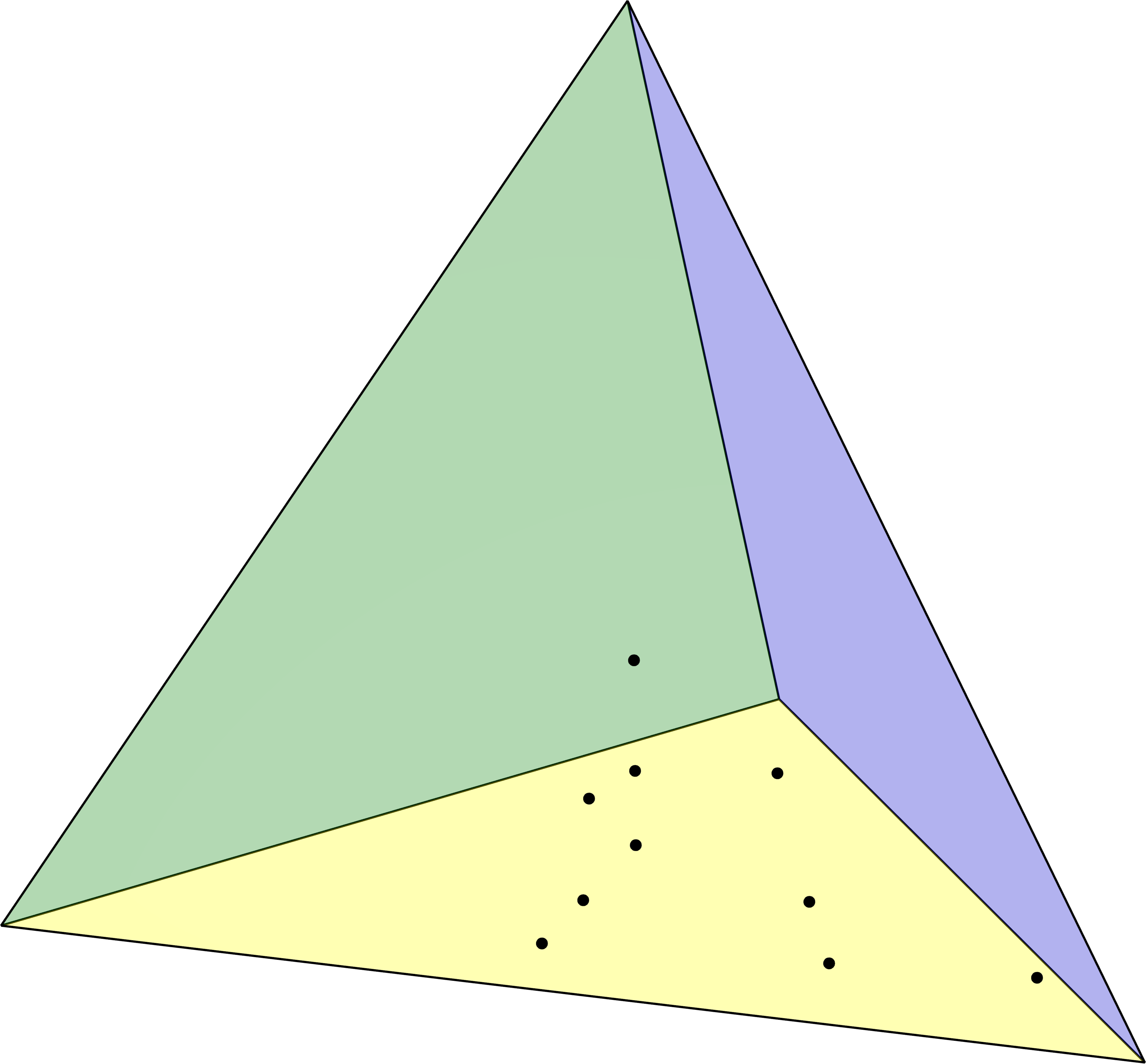}
		\caption{\label{fig:unique sym tet} Unique symmetry orbits}
	\end{subfigure}\hfill
	\begin{subfigure}{0.30\textwidth}
		\centering
		\includegraphics[scale=0.0130]{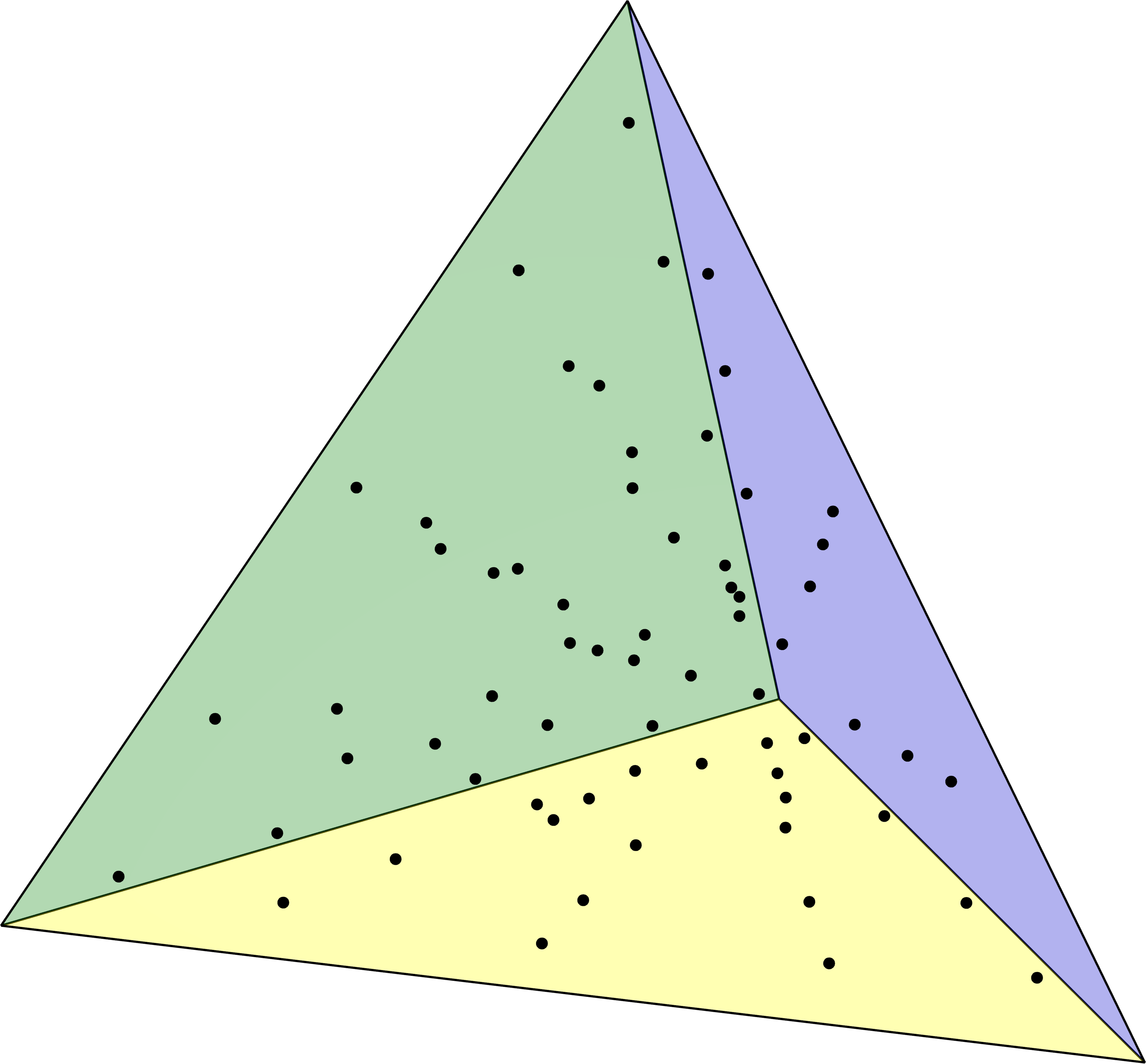}
		\caption{\label{fig:init tet} Initial node distribution}
	\end{subfigure}
	\caption{\label{fig:initial guess tri tet} Generation of initial nodes for the degree $q=4$ line-LG quadrature rule on triangles and tetrahedra.}
\end{figure}

The initial guesses for the PI quadrature rules are generated by forming tensor-product structures on quadrilateral/hexahedral subdomains of the simplices using the LG nodes lying in the first half of the line reference element. The construction procedure is best explained through an example. Consider the construction of the degree $q=8$ line-LG quadrature rule on the triangle; we first map the nodes of the degree $q=9$ one-dimensional LG quadrature rule that lie in the first half of the line reference element onto the line connecting the bottom left vertex to the bottom edge midpoint. We apply the same mapping to the line connecting the left edge midpoint to the centroid of the triangle.  Then we connect the corresponding LG nodes on the two lines and, once again, apply the mapping to each of the connecting lines to find the initial nodes on a quadrilateral subdomain of the triangle, see \cref{fig:lg line tri}. The unique symmetry orbits can be found using the coordinates of the nodes that lie in one of the triangles obtained by dividing one of the three quadrilateral subdomains into two triangles, as shown in \cref{fig:unique sym tri}. Denoting the number of nodes in the unique orbits and their nodal coordinate matrix by $n_q^{\prime}$ and $\bm{y}^{\prime}\in \IR{d \times n_q^{\prime}}$, respectively, we can compute the barycentric coordinates using the relation in \cref{eq:x=T lambda} and the fact that the sum of the barycentric coordinates of each node is unity, \ie,
\begin{equation}
	\bm{\lambda} = (\ubar{\T}^{T})^{-1}\ubar{\bm{y}}^\prime,
\end{equation}
where $(\ubar{\cdot})$ denotes concatenation of an additional row vector of ones at the bottom of a matrix. The complete nodal distribution on the triangle shown in \cref{fig:init tri} can be computed using the barycentric coordinates and the symmetry relations in the first line of \cref{eq:sym tri}. For the tetrahedron, the unique symmetry orbits are provided by the nodes that lie in one of the six tetrahedra obtained by dividing one of the four hexahedral subdomains of the reference tetrahedron into six tetrahedra, see \cref{fig:unique sym tet}. 

The initial weight parameters associated with the unique symmetry orbits are set to a constant value of $(2/q)^3$ in all cases, as the solution to the line-LG quadrature problem is observed to be relatively less sensitive to the values of the weights than the location of the nodes. 

\section{The Line-LG Quadrature Rule}\label{sec:Line-LG}
With the node parameters and weights set for each symmetry group, we proceed to solve the minimization problem in \cref{eq:min prob} using the LM algorithm. We derived line-LG quadrature rules up to degree $q=84$ and $q=40$ on the triangle and tetrahedron, respectively. The reason we stopped at $q=84$ for the triangle is because the outputs of the classic Gamma function used in the construction of the PKD basis functions are out of the range of numbers that can be represented in double precision. For the tetrahedron, memory limitations become significant at quadrature degrees slightly larger than forty. While both of these limitations can be circumvented, we chose not to pursue this task and have left it for future work.

The type and number of symmetry orbits in the line-LG quadrature rules can be determined by studying the structure of the initial nodal distribution. The one-dimensional LG quadrature rule is exact for polynomials up to degree $q=2n_1 - 1$, where $n_1$ is the number of quadrature points on the reference line. Defining $m= (n_1 \bmod 2)$ and $n_r = (n_1- m)/2$, the number of orbits in each type of symmetry group of the triangle can be calculated as
\begin{align}
	|S_{1}|= m, \quad |S_{21}| = (1 +m) n_r,   \quad  |S_{111}| = (n_r^2 - n_r)/2.
\end{align} 
Similarly, the type and numbers of symmetry orbits on the reference tetrahedron are 
\begin{equation}
	\begin{aligned}
		|S_{1}| =& m, \qquad |S_{31}| =  (1 + m) n_r, \qquad |S_{22}| = m n_r,\\
		|S_{211}| =& \frac{1+2m}{1+m} (n_r^2 - n_r),  \qquad |S_{1111}| = \frac{(n_r - 1)^3 - n_r +1}{6}.
	\end{aligned}
\end{equation}
Note that the above numbers of orbits are functions of $n_1$ only; hence, they provide information for the initial guesses generated using LG quadrature rules with $n_1$ nodes. The actual derivation of a degree $q$ line-LG quadrature rule may use $n_1$ values that correspond to a higher degree one-dimensional LG quadrature rule, as explained in the next section.

\subsection{Efficiency of the Line-LG Quadrature Rule}\label{sec:efficiency line-LG}
The efficiency of the line-LG quadrature rules depends on the number of one-dimensional quadrature points, $n_1$, used in their construction and the polynomial degree they can integrate exactly. On the triangle, we have observed some quadrature accuracy patterns associated with the degree and values of $n_1$. If $q-1$ is divisible by $4$, then only $n_1 = \lfloor q/2\rfloor + 1$ nodes are needed to derive the line-LG quadrature rule. This case coincides with the presence of a mid-point in the one-dimensional LG quadrature rule.  If $q$ is even or $q\ge30$, then $n_1=\lfloor q/2\rfloor + 1$ nodes are required to derive the line-LG quadrature rules. Finally, if $q<30$, odd, and $q-1$ is not divisible by $4$, then $n_1=\lfloor q/2 \rfloor+2$ nodes are required.  For example, the degree $8$, $9$, $10$, and $11$ line-LG quadrature rules require initial guesses constructed with $n_1$ values of $5$, $5$, $6$, and $7$, respectively. This has an impact on the efficiency of the rules, in the sense of definition \cref{eq:eff quad}, as, for example, the $q=9$ rule has larger lower bounds on the number of nodes, but it has the same number of nodes as the $q=8$ rule. 

On the tetrahedron, a more consistent pattern is observed. For all values of $q$ except $3$, $7$, and $11$, we were able to find line-LG quadrature rule using $n_1=\lfloor q/2 \rfloor + 1$ nodes to construct the initial guesses. For the exceptions, $n_1=\lfloor q/2 \rfloor + 2$ nodes are required. The solutions for the line-LG quadrature rules are obtained in a few iterations of the LM algorithm. For example, the $q=20$ line-LG rule on the tetrahedron requires less than $10$ iterations.  

\begin{figure}[t]
	\centering
	\begin{subfigure}[c]{0.03\textwidth}
		\centering
		\caption{\label{fig:eff lg line tri} }
	\end{subfigure}
	\begin{minipage}[c]{0.95\textwidth}
		\includegraphics[scale=0.485]{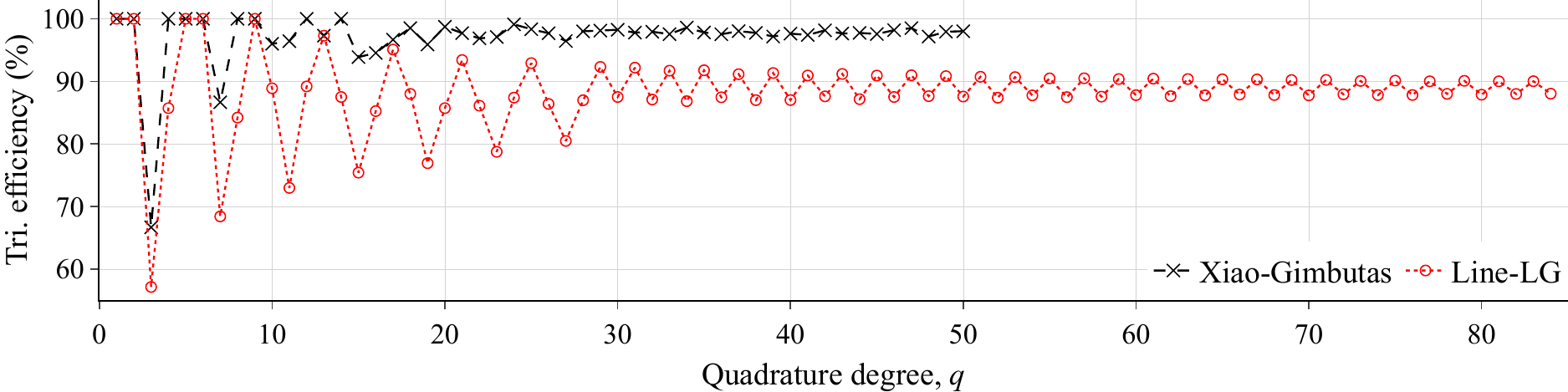}
	\end{minipage}
	\\
	\begin{subfigure}[c]{0.03\textwidth}
		\centering
		\caption{\label{fig:eff lg line tet} }
	\end{subfigure}
	\begin{minipage}[c]{0.95\textwidth}
		\includegraphics[scale=0.485]{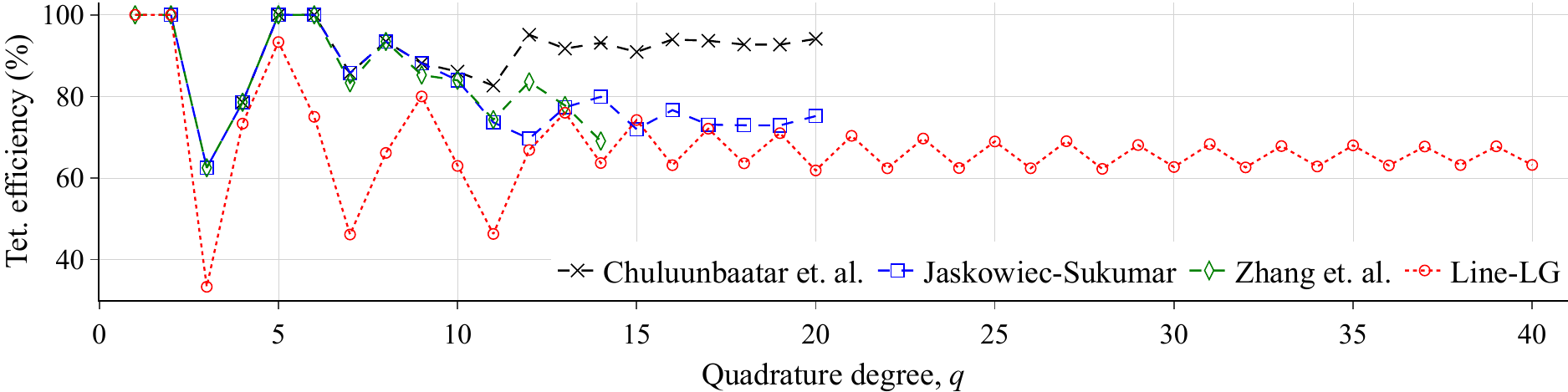}
	\end{minipage}
	
	\caption{\label{fig:eff lg line} Efficiency of the line-LG quadrature rules relative to existing symmetric quadrature rules on triangles (a) and tetrahedra (b).}
\end{figure}

Despite their tensor-product-like initial guess structures and ease of derivation, the line-LG quadrature rules are surprisingly efficient. \cref{fig:eff lg line} shows their efficiency relative to other known quadrature rules in the literature. For triangles the line-LG rules are above $80\%$ efficient for most degrees, and are close to $90\%$ efficient at very high degrees. For the tetrahedron, they are more than $60\%$ efficient for most degrees, which is significant given that many of the existing rules in the literature, except the recent additions of Chuluunbaatar \etal \cite{chuluunbaatar2022pi}, are less than $80\%$ efficient for degrees greater than ten. 

\section{Node Elimination}\label{sec:node elim}
Node elimination is a technique used to create quadrature rules with fewer nodes, starting from existing known rules \cite{xiao2010numerical}. The idea is straightforward; given a quadrature rule, eliminate an orbit and attempt to solve the minimization problem in \cref{eq:min prob} to machine precision. If the attempt is unsuccessful, reintroduce the eliminated orbit. If successful, a new quadrature rule is obtained, and the process is repeated until no further orbit can be eliminated. When there are several orbits that can be eliminated at a given stage, different criteria can be used to decide which orbit to eliminate. For example, Xiao and Gimbutas \cite{xiao2010numerical} proposed a significance index for each node, which is the sum or weighted sum of the squares of the basis functions evaluated at each node. The node with the smallest significance index is considered to have the lowest contribution to the integral approximation; hence, it is eliminated first. Another approach is to eliminate the node with the smallest weight, which is appealing from the perspective of collocated high-order methods since this will reduce the norm of the inverse of the mass or norm matrix and, consequently, the conditioning of the system matrices resulting from discretizations of PDEs. One can also eliminate orbits such that the minimum distance between nodes is increased, which can improve the CFL number for numerical solutions of PDEs using explicit time-marching methods. It may also be desirable to eliminate nodes that are very close to the element boundaries. While all of these options appear reasonable, it is not clear a priori which approach will ultimately result in the minimum number of nodes. The decision at each stage will affect the set of orbits that can be eliminated in subsequent stages, thereby determining the total number of nodes that can be attained. Due to the combinatorial nature of the elimination process, it is impractical to test all possible paths, particularly at higher degrees. In this work, we have not made attempts to investigate the optimal node elimination strategy. However, based on limited experiments, we have opted to eliminate orbits that produce nodes closest to the facets or that have the smallest weights. For triangles, we have used both of these approaches and neither is conclusively better than the other. For the tetrahedron, the first approach is used in most cases, but when efficiency of a derived quadrature is low, both approaches are applied to check if improvements can be made. 

The node elimination process takes into account the structure of the lower-bound estimates presented in \cref{eq:tri bound} and \cref{eq:tet bound}. In particular, for the first two outer elimination iterations, see \cref{fig:flowchart}, the number of nodes in each symmetry orbit is kept equal to or greater than the number of nodes presented in the estimates. This increases the possibility of attaining quadrature rules that match the estimated bounds or have a few more orbits than the estimates. Furthermore, at each stage in the first two outer iterations, the type of orbit that supports the highest number of nodes is eliminated until no further orbit of that kind can be eliminated. For example, at any stage in the first two outer iterations, the $S_{111}$ orbits are ranked based on proximity to the facets or based on their weights, and the highest ranking orbit is eliminated before any of the $S_{21}$ orbits are considered. After the second outer iteration, all types of symmetry orbits are ranked at once and the highest ranking orbit is eliminated. The node elimination process is sensitive to the value of $\nu$ used in the LM algorithm, which is defined in \cref{eq:step}. Hence, values of $\nu$ in orders of $10$ ranging from $10^{-8}$ to $10^4$ are used. Admittedly, the process of node elimination requires changing the above parameters manually, starting from higher degree line-LG rules, or, in rare cases, combining two or more line-LG quadrature rules to improve efficiency. Automation and robust elimination algorithms are lacking and are worth pursuing in order to further improve efficiency and extend the proposed approach to dimensions greater than three. 

\section{New Symmetric PI Quadrature Rules}\label{sec:new quad}
\begin{figure}[!t]
	\centering
	\begin{subfigure}{0.24\textwidth}
		\centering
		\includegraphics[scale=0.26]{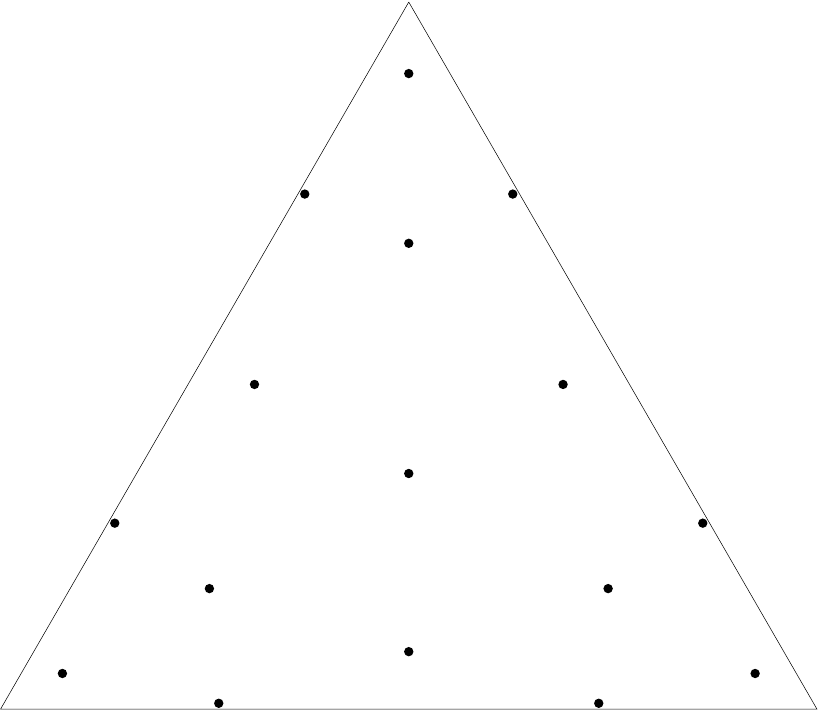}
		\caption{$q=8$, $n_q=16$}
	\end{subfigure}\hfill
	\begin{subfigure}{0.24\textwidth}
		\centering
		\includegraphics[scale=0.26]{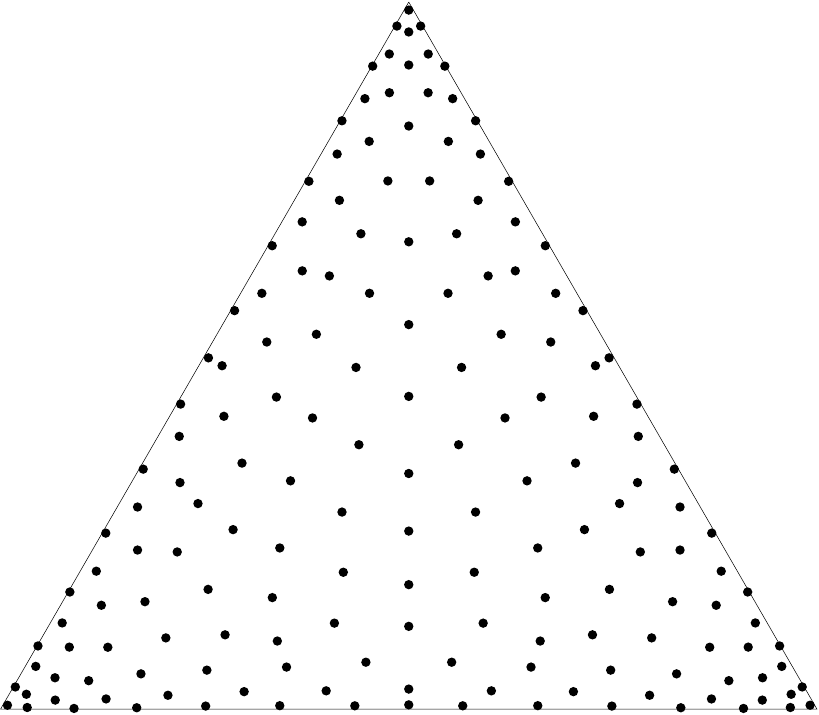}
		\caption{$q=30$, $n_q=169$}
	\end{subfigure}\hfill
	\begin{subfigure}{0.24\textwidth}
		\centering
		\includegraphics[scale=0.26]{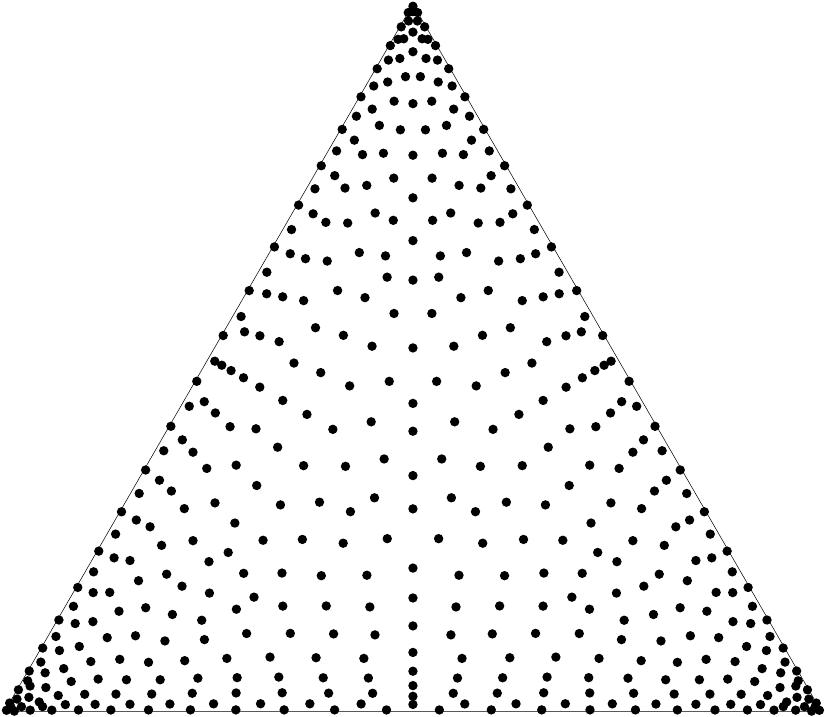}
		\caption{$q=50$, $n_q=448$}
	\end{subfigure}\hfill
	\begin{subfigure}{0.24\textwidth}
		\centering
		\includegraphics[scale=0.26]{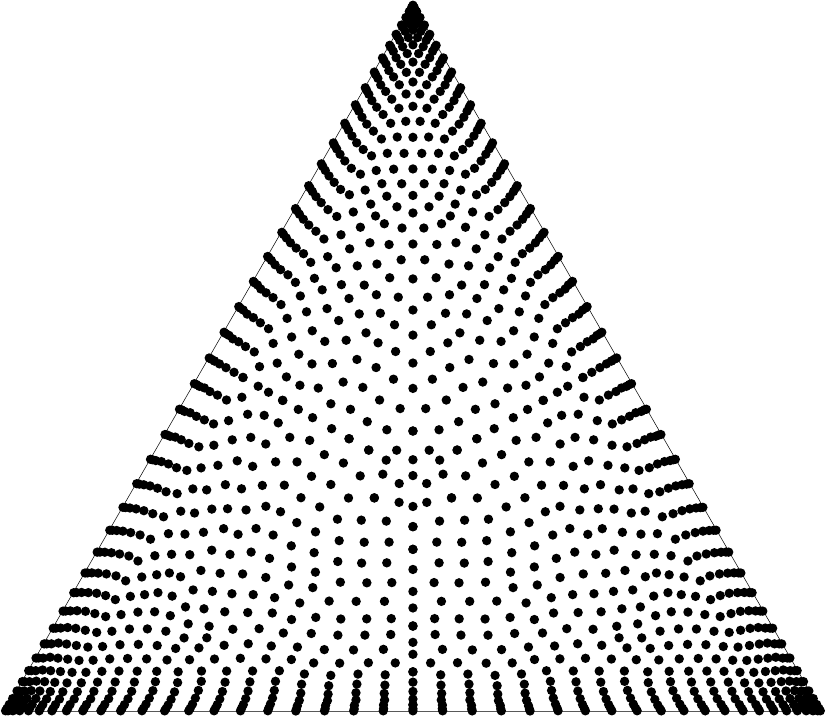}
		\caption{$q=84$, $n_q=1261$}
	\end{subfigure}
	\\ \vspace{0.2cm}
	\begin{subfigure}{0.24\textwidth}
		\centering
		\includegraphics[scale=0.0425]{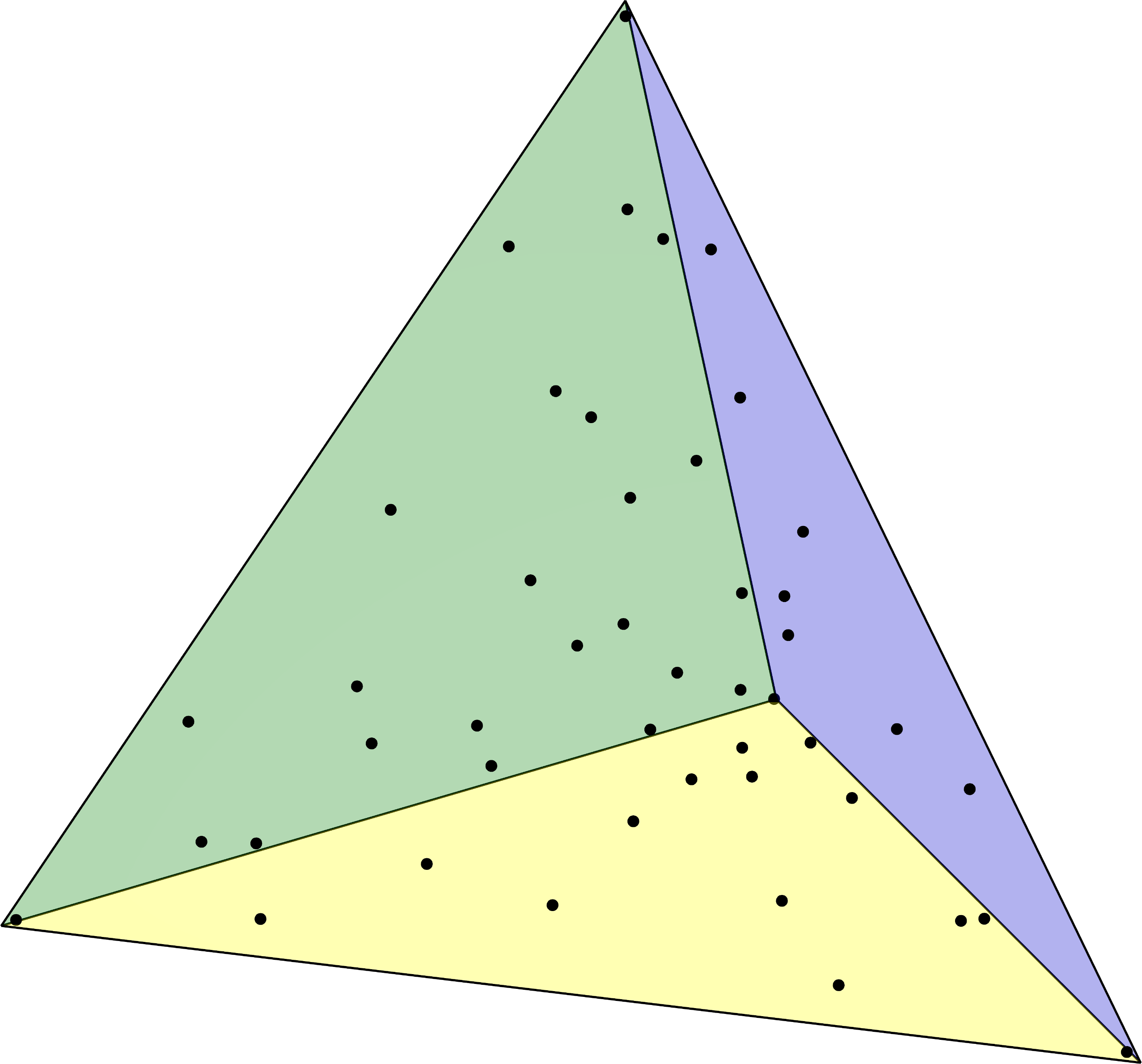}
		\caption{$q=8$, $n_q=46$}
	\end{subfigure}\hfill
	\begin{subfigure}{0.24\textwidth}
		\centering
		\includegraphics[scale=0.0425]{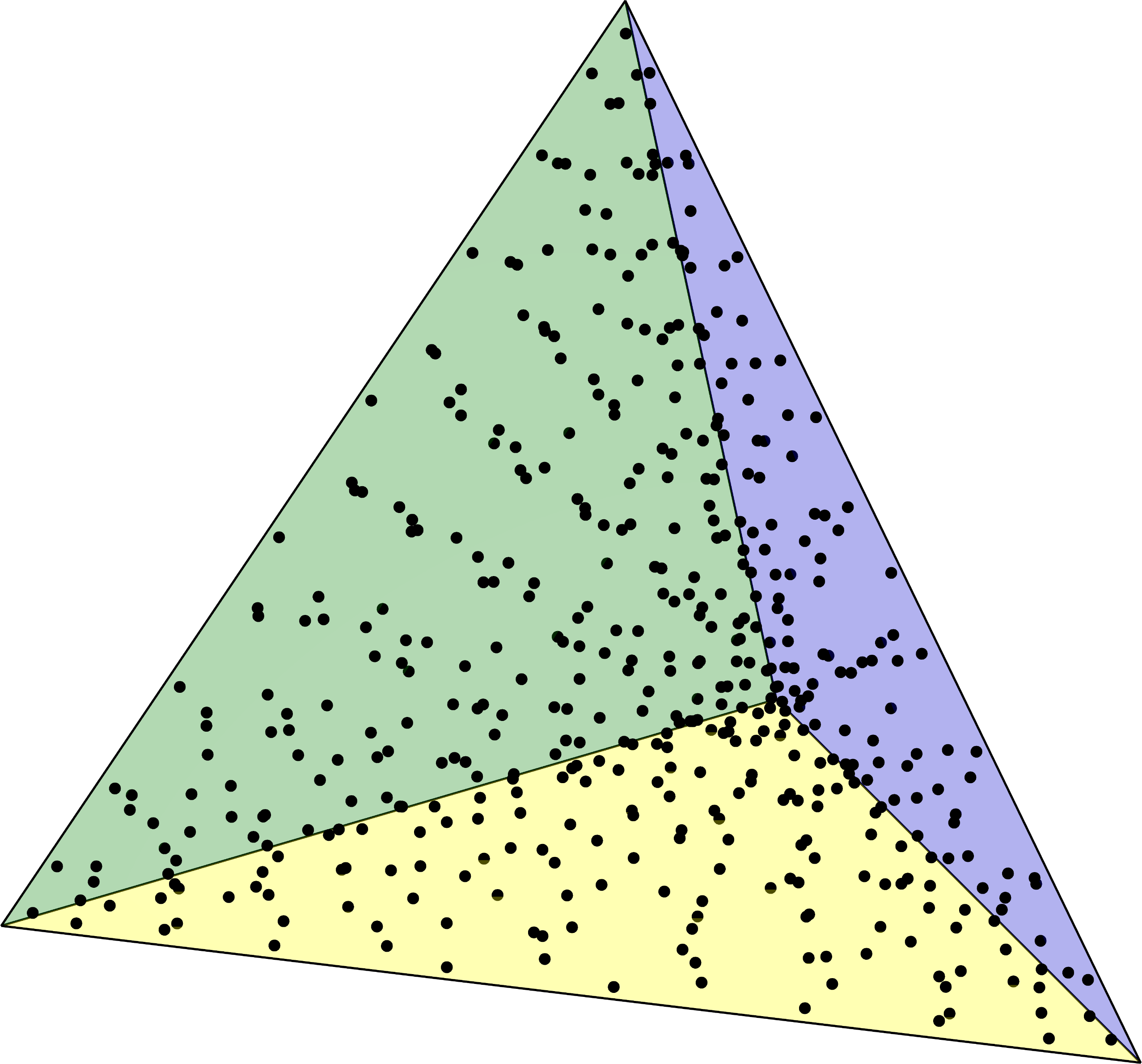}
		\caption{$q=20$, $n_q=469$}
	\end{subfigure}\hfill
	\begin{subfigure}{0.24\textwidth}
		\centering
		\includegraphics[scale=0.0425]{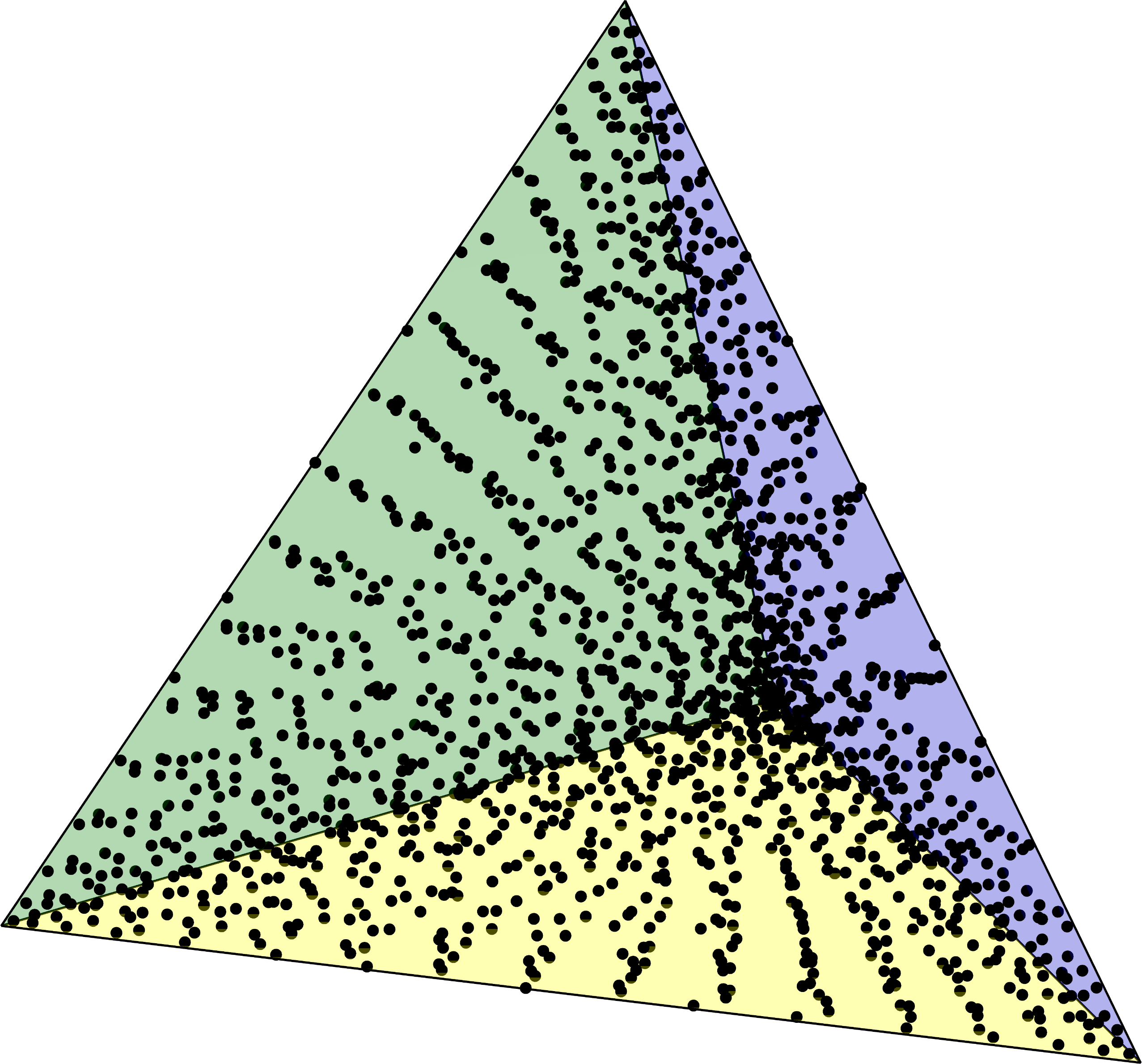}
		\caption{$q=31$, $n_q=1667$}
	\end{subfigure}\hfill
	\begin{subfigure}{0.24\textwidth}
		\centering
		\includegraphics[scale=0.0425]{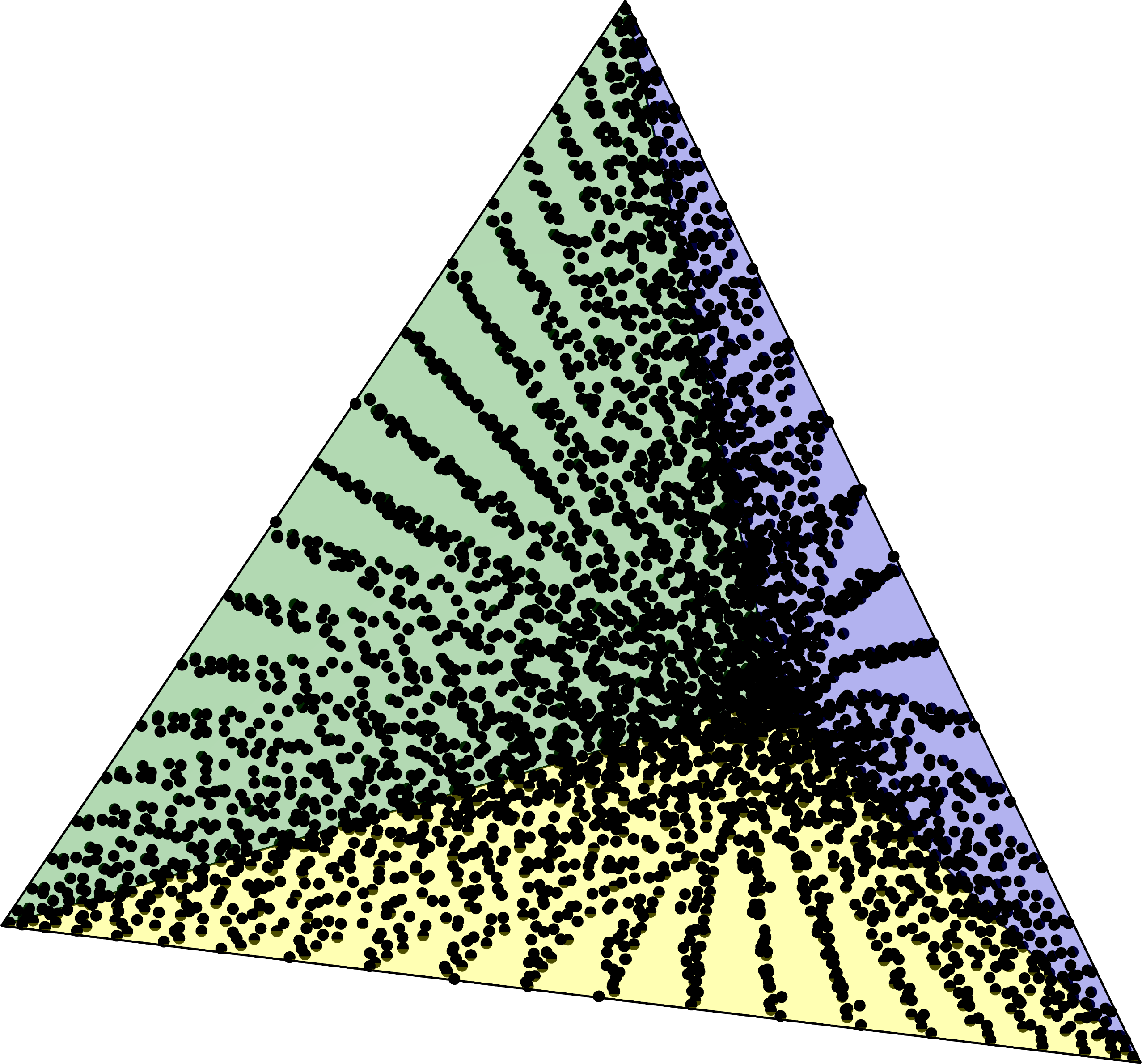}
		\caption{$q=40$, $n_q=3815$}
	\end{subfigure}
	\caption{\label{fig:tri tet nodes elim} Examples of symmetric PI quadrature rules on triangles and tetrahedra.}
\end{figure}

By applying the strategies outlined in the previous sections, we derived more efficient quadrature rules than the line-LG rules on both triangles and tetrahedra. The number of nodes for the quadrature rules on triangles are shown in \cref{tab:numnodes tri}. While a few of the newly derived rules on the triangle have more nodes than the existing rules, most of them have either an equal number or fewer nodes. \cref{tab:numnodes tet} shows that many of the derived rules on the tetrahedron have fewer nodes than those presented by Ja{\'s}kowiec and Sukumar \cite{jaskowiec2021high}, but all of the rules up to degree twenty have an equal number or more nodes than the recent results of Chuluunbaatar \etal  \cite{chuluunbaatar2022pi}. To the authors' knowledge all of the rules greater than degree twenty are new. The nodal locations of some of the derived quadrature rules are shown in \cref{fig:tri tet nodes elim}.
\begin{table*}[!t]
	\footnotesize
	\centering
	\caption{\label{tab:numnodes tri} Number of nodes of the newly derived symmetric PI quadrature rules on triangles compared to existing rules. New quadrature rules with fewer nodes than existing rules and those for higher degrees than previously available are underlined.}
	\begin{threeparttable}
		\setlength{\tabcolsep}{0.1em}
		\renewcommand*{\arraystretch}{1.2}
		\begin{tabular}{ccc@{\hspace{1.0em}}ccc@{\hspace{1.0em}}ccc@{\hspace{1.0em}}ccc@{\hspace{1.0em}}ccc@{\hspace{1.0em}}cc@{\hspace{1.0em}}cc@{\hspace{1.0em}}cc@{\hspace{1.0em}}cc }
			\toprule
			$q$ & \cite{xiao2010numerical} & new &  $q$ & \cite{xiao2010numerical} & new &  $q$ & \cite{xiao2010numerical} & new &  $q$ & \cite{xiao2010numerical} & new &  $q$ & \cite{xiao2010numerical} & new &  $q$ & new &  $q$ & new & $q$ & new & $q$ & new \\
			\midrule
			1  &  1  &  1  & 11 & 28 & 28 & 21 & 87 & 87 & 31 & 181 & 181 & 41 & 309 & \underline{304} & 51 & \underline{468} & 61 & \underline{670} & 71 & \underline{913} & 81 & \underline{1174} \\
			2  &  3  &  3  & 12 & 33 & 33 & 22 & 96 & 97 & 32 & 193 & \underline{192} & 42 & 324 & \underline{321} & 52 & \underline{484} & 62 & \underline{693} & 72 & \underline{931} & 82 & \underline{1188} \\
			3  &  6  &  6  & 13 & 37 & 37 & 23 & 103 & 103 & 33 & 204 & \underline{202} & 43 & 339 & \underline{337} & 53 & \underline{504} & 63 & \underline{709} & 73 & \underline{936} & 83 & \underline{1222} \\
			4  &  6  &  6  & 14 & 42 & 42 & 24 & 112 &\underline{111} & 34 & 214 & 214 & 44 & 354 & \underline{349} & 54 & \underline{522} & 64 & \underline{724} & 74 & \underline{1000} & 84 & \underline{1261} \\
			5  &  7  &  7  & 15 & 49 & 49 & 25 & 120 & 121 & 35 & 228 & \underline{226} & 45 & 370 & \underline{367} & 55 & \underline{541} & 65 & \underline{748} & 75 & \underline{1024} &  &  \\
			6  & 12 & 12 & 16 & 55 & 55 & 26 & 130 & 130 & 36 & 243 & \underline{240} & 46 & 385 & \underline{382} & 56 & \underline{561} & 66 & \underline{793} & 76 & \underline{1033} &  & \\
			7  & 15 & 15 & 17 & 60 & 60 & 27 & 141 & \underline{139} & 37 & 252 & \underline{250} & 47 & 399 & \underline{397} & 57 & \underline{579} & 67 & \underline{808} & 77 & \underline{1081} &  &  \\
			8  & 16 & 16 & 18 & 67 & 67 & 28 & 150 & \underline{148} & 38 & 267 & \underline{265} & 48 & 423 & \underline{415} & 58 & \underline{598} & 68 & \underline{829} & 78 & \underline{1099} &  &  \\
			9  & 19 & 19 & 19 & 73 & 73 & 29 & 159 & 159 & 39 & 282 & \underline{277} & 49 & 435 & \underline{430} & 59 & \underline{628} & 69 & \underline{865} & 79 & \underline{1117} &  &  \\
			10 & 25 & 25 & 20 & 79 & 79 & 30 & 171 & \underline{169} & 40 & 295 & \underline{292} & 50 & 453 & \underline{448} & 60 & \underline{649} & 70 & \underline{868} & 80 & \underline{1129} &  &  \\
			\bottomrule
		\end{tabular}
	\end{threeparttable}
\end{table*}
\begin{table*}[!t]
	\footnotesize
	\centering
	\caption{\label{tab:numnodes tet} Number of nodes of the newly derived symmetric PI quadrature rules on tetrahedra compared to existing rules. New quadrature rules for higher degrees than previously available are underlined.}
	\begin{threeparttable}
		\setlength{\tabcolsep}{0.1em}
		\renewcommand*{\arraystretch}{1.2}
		\begin{tabular}{cccc@{\hspace{1.0em}}cccc@{\hspace{1.0em}}cccc@{\hspace{1.0em}}cccc@{\hspace{1.0em}}cc@{\hspace{1.0em}}cc@{\hspace{1.0em}}cc@{\hspace{1.0em}}cc}
			\toprule
			$q$ & \cite{jaskowiec2021high} & \cite{chuluunbaatar2022pi}& new &  $q$ & \cite{jaskowiec2021high} & \cite{chuluunbaatar2022pi}& new & $q$ & \cite{jaskowiec2021high} & \cite{chuluunbaatar2022pi}& new & $q$ & \cite{jaskowiec2021high} & \cite{chuluunbaatar2022pi}& new &  $q$ & new  & $q$ & new  & $q$& new & $q$ & new \\
			\midrule
			1  &  1  	&   	& 1 & 	6  &  24 & 24 & 24&  11 & 110 &98& 100 & 		 16 & 304& 248 & 259& 			21 & \underline{540} & 26 &  \underline{992} & 31 &  \underline{1667}  		& 36 &  \underline{2652}  \\
			2  &  4  	&    	&4 &  7  & 35 & 35& 35&  	12 & 168 & 123 &124 & 		17 &  364& 284& 299& 		22 &  \underline{616} & 27 &  \underline{1148}  &32 &  \underline{1768}  		& 37 &  \underline{2671} \\
			3  &  8  	&   	& 8& 8  & 46 & 46 & 46& 	13 & 172& 145 &145 &		18 & 436& 343& 370&      	23 &  \underline{706} &28 &  \underline{1200} &  33 & \underline{1913}  		&38 &  \underline{3296} \\
			4  &  14  &  14  &14 & 9  & 59 & 59 & 61&	14 & 204 & 175 &181& 		19 & 487& 383& 408& 		24 &  \underline{756} &  29 &  \underline{1405} & 34 &  \underline{2175}   		&39 &  \underline{3436}\\
			5  &  14  &  14 & 14& 10 & 81 & 79 & 79&	15 & 264 & 209 &216 & 		20 & 552& 441 & 469& 		25 &  \underline{904} &30 &  \underline{1564} &  35 &  \underline{2352} 		&40 &  \underline{3815}  \\
			\bottomrule
		\end{tabular}
	\end{threeparttable}
\end{table*}

\begin{figure}[t!]
	\centering
	\begin{subfigure}[c]{0.03\textwidth}
		\centering
		\caption{\label{fig:eff lg tri} }
	\end{subfigure}
	\begin{minipage}[c]{0.95\textwidth}
		\includegraphics[scale=0.485]{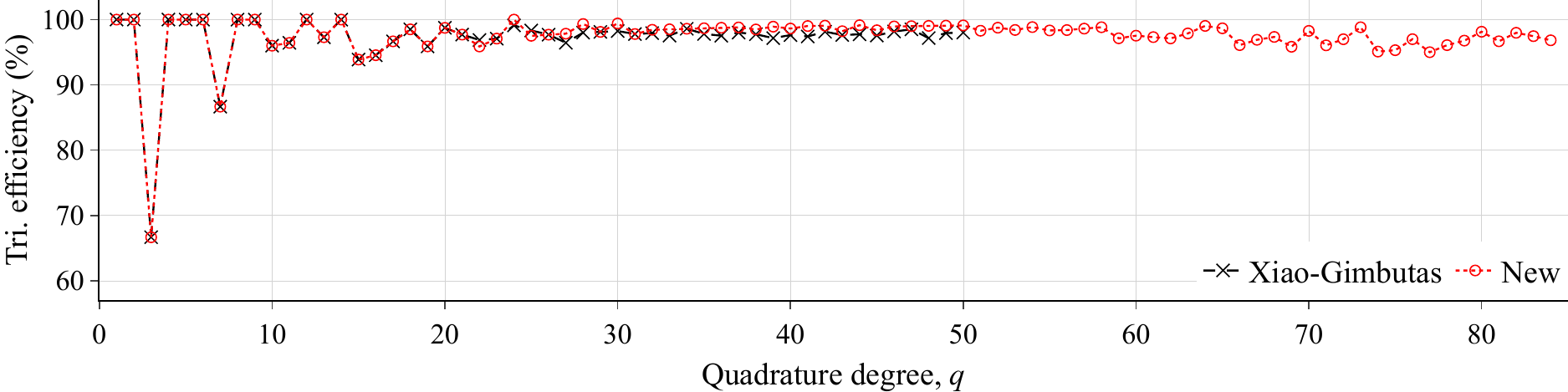}
	\end{minipage}
	\\
	\vspace{0.2cm}
	\begin{subfigure}[c]{0.03\textwidth}
		\centering
		\caption{\label{fig:eff lg tet} }
	\end{subfigure}
	\begin{minipage}[c]{0.95\textwidth}
		\includegraphics[scale=0.485]{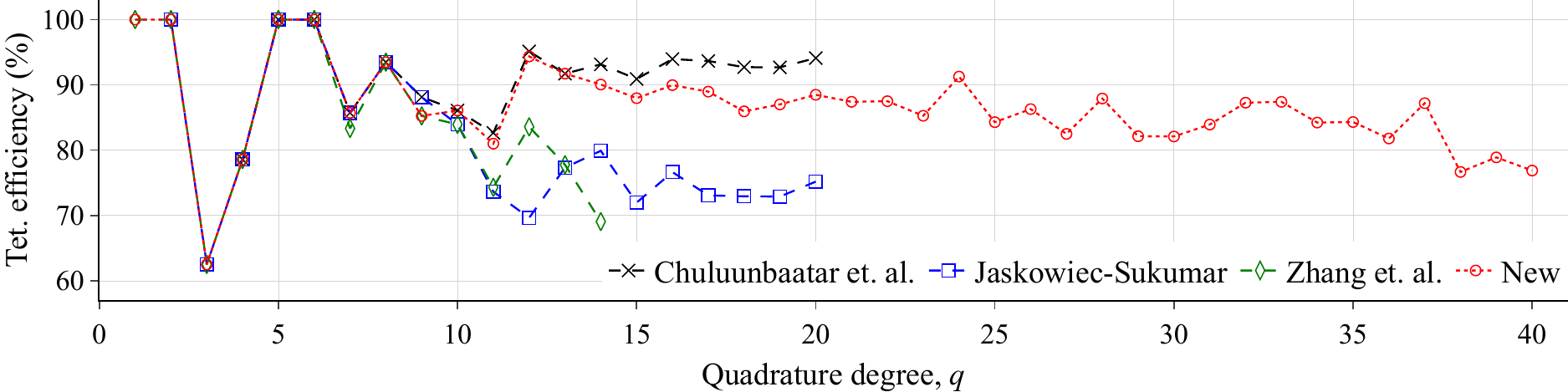}
	\end{minipage}
	
	\caption{\label{fig:eff line} Efficiency of the quadrature rules derived in this work relative to existing symmetric quadrature rules on triangles (a) and tetrahedra (b).}
\end{figure}

The efficiency of the quadrature rules obtained after applying node elimination is shown in \cref{fig:eff line} for both the triangle and tetrahedron. Most of the quadrature rules for the triangle are more than $95\%$ efficient, while all but a few are more than $80\%$ efficient for the tetrahedron. It is evident from \cref{fig:eff lg tet} that most of the newly derived quadrature rules on the tetrahedron have comparable efficiency relative to the most efficient existing rules, with at most a $7\%$ reduction in efficiency.

All of the quadrature rules are provided in the supplementary material accompanying this paper. The line-LG quadrature rules are also included, which could play an important role in future developments of robust node elimination algorithms and derivation of more efficient quadrature rules.

\section{Numerical Results}\label{sec:numerical}
We apply the derived quadrature rules to compute the following integrals of highly oscillatory functions,
\begin{align}
\fn{I}_{2} &= \int_{0}^{1}\int_{0}^{1}\sin(48\pi x_1^8) \cos(48 \pi x_2^5)\dd{x}\dd{y},\label{eq:integ I}\\
\fn{I}_{3} &= \int_{0}^{1}\int_{0}^{1}\int_{0}^{1}\sin(16\pi x_1^8) \cos(16 \pi x_2^5) \cos(16 \pi x_3^7)\dd{x}\dd{y}\dd{z},\label{eq:integ J}
\end{align}
in two and three dimensions, respectively. The exact values of these integrals can be obtained using symbolic mathematics libraries, and the first sixteen significant digits are $0.03116210698718051$ and $0.02288392144769807$, respectively. To test the accuracy of the quadrature rules, we constructed meshes for the unit square and unit cube that maintained an approximately equal number of degrees of freedom across different quadrature degrees.  For the square, the meshes were constructed such that the total number of nodes was within $7\%$ of $567\,450$.  For the cube, the tetrahedral meshes were constructed such that the total number of nodes was within $7\%$ of $158\,779\,350$. The stated reference total number of nodes correspond to having $450$ degree $q=84$ triangular elements and $41,154$ degree $q=40$ tetrahedral elements with the new quadrature rules (we remind the reader that the label ``new'' refers to the quadrature rules obtained after node elimination).
\begin{figure}[t!]
	\centering
	\begin{subfigure}[c]{0.03\textwidth}
		\centering
		\caption{\label{fig:error tri} }
	\end{subfigure}
	\begin{minipage}[c]{0.95\textwidth}
		\includegraphics[scale=0.485]{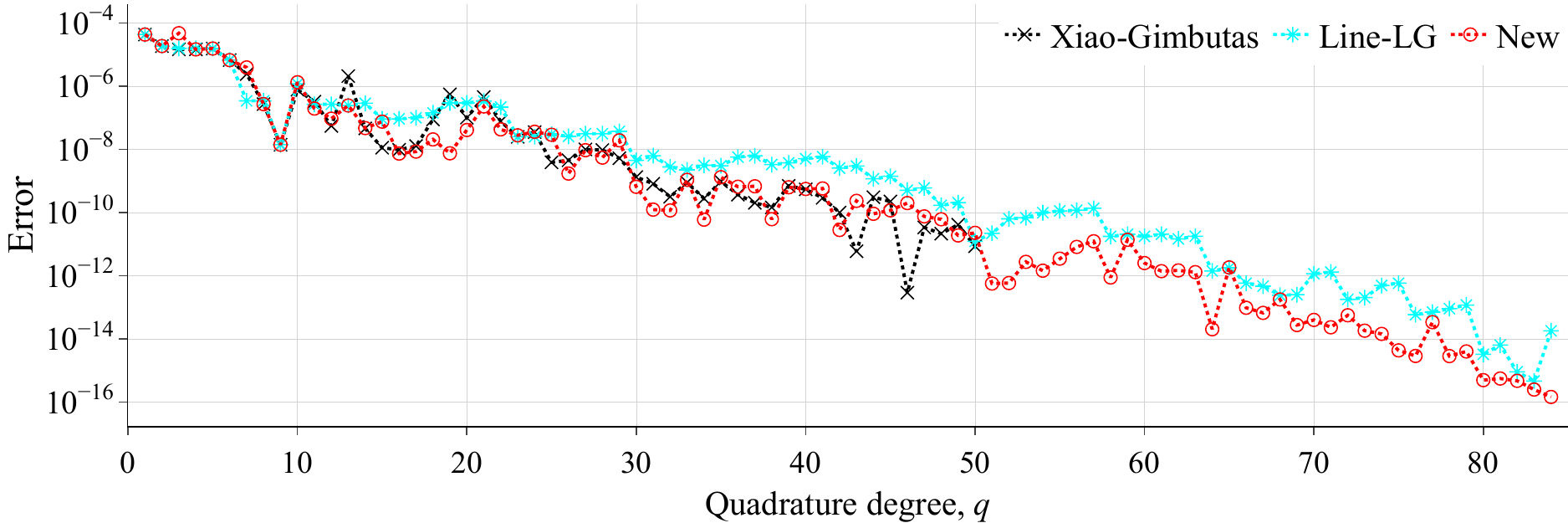}
	\end{minipage}
	\\
	\vspace{0.2cm}
	\begin{subfigure}[c]{0.03\textwidth}
		\centering
		\caption{\label{fig:error tet} }
	\end{subfigure}
	\begin{minipage}[c]{0.95\textwidth}
		\includegraphics[scale=0.485]{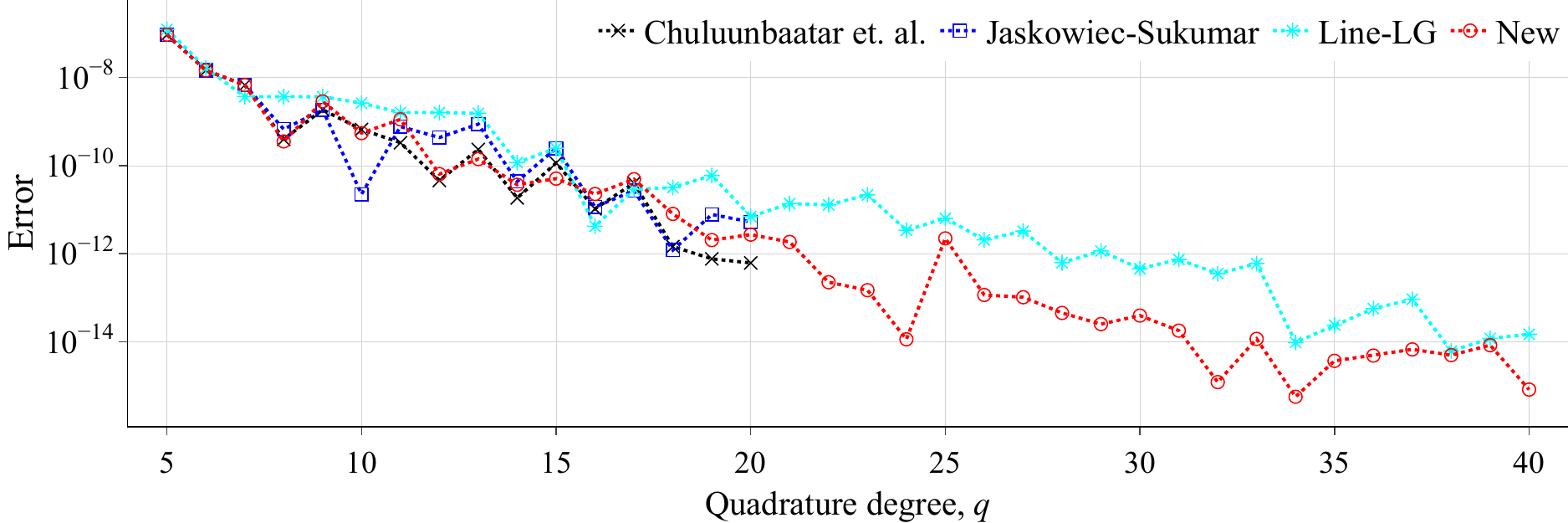}
	\end{minipage}
	
	\caption{\label{fig:error} Quadrature approximation absolute errors on triangles (a) and tetrahedra (b). The meshes used for each quadrature rule have approximately equal number of degrees of freedom.}
\end{figure}

\cref{fig:error} shows the absolute values of the quadrature approximation errors of the integrals in \cref{eq:integ I} and \cref{eq:integ J}. The new quadrature rules yield better approximations than the line-LG rules in most cases. Furthermore, they show comparable accuracy relative to the existing rules of \cite{xiao2010numerical} and \cite{jaskowiec2021high} in all but some isolated cases. In particular, the new quadrature rules yield significantly less accurate results for the degree $10$ case in three dimensions and degree $43$ and $46$ in two dimensions. 

To investigate the three-dimensional case further, we apply the degree 8 through 12 quadrature rules to approximate the integral of a less oscillatory function,
\begin{equation}
\fn{J}_{3} = \int_{0}^{1}\int_{0}^{1}\int_{0}^{1}\sin(3\pi x)\sin(5\pi y)\sin(3\pi z) \dd{x}\dd{y}\dd{z}= \frac{8}{45\pi^3}. \label{eq:integ smooth}
\end{equation}
\begin{table*} [!t]
    \footnotesize
    \caption{\label{tab:grid conv} Grid convergence study of quadrature approximation errors on tetrahedral meshes.}
    \centering
    \setlength{\tabcolsep}{0.3em}
    \renewcommand*{\arraystretch}{1.2}
    \begin{tabular}{ll@{\hspace{1em}}ll@{\hspace{1.75em}}ll@{\hspace{1.75em}}ll@{\hspace{1.75em}}ll}
        \toprule
        $q$ & \multirow{2}{*}{$ \#elem $ } & \multicolumn{2}{c}{\cite{chuluunbaatar2022pi}} & \multicolumn{2}{c}{\cite{jaskowiec2021high}} &  \multicolumn{2}{c}{Line-LG} & \multicolumn{2}{c}{New} \\
        \cmidrule(l{0em}r{1.5em}){3-4} \cmidrule(l{0em}r{1.5em}){5-6} \cmidrule(l{0em}r{1.5em}){7-8} \cmidrule(l{0em}r{0.25em}){9-10}
        & & {error}  &{rate} & {error}  & {rate} & {error} & {rate} & {error} & {rate}\\
        \midrule
        \multirow[t]{4}{*}{$ 8 $}&    $ 6^3\times 6 $ 		& $5.5338e-11$&  &$7.9018e-10$&  & $1.3511e-09$&  & $1.8845e-09$&  \\  
        &    $ 7^3\times 6 $ 		& $1.2141e-11$& $9.84$ & $1.5439e-10$& $10.59$ & $2.6181e-10$& $10.65$ & $3.6643e-10$& $10.62$\\  
        &    $ 8^3\times 6 $		& $3.2247e-12$& $9.93$ & $3.8351e-11$& $10.43$ & $6.4692e-11$& $10.46$ & $9.0741e-11$& $10.45$ \\   
		&    $ 9^3\times 6 $		& $9.9720e-13$& $9.96$ & $1.1363e-11$& $10.33$ & $1.9099e-11$& $10.36$ & $2.6830e-11$& $10.35$\\    
        \addlinespace
		\multirow[t]{4}{*}{$ 9 $}&    $ 6^3\times 6 $ 		& $1.2350e-09$ &  & $1.2339e-09$&  &  $1.3493e-09$&  & $1.3491e-09$&   \\  
        &    $ 7^3\times 6 $ 		& $2.3928e-10$& $10.65$ & $2.3914e-10$& $10.64$ &  $2.6146e-10$& $10.64$ & $2.6136e-10$& $10.65$ \\  
        &    $ 8^3\times 6 $		& $5.9118e-11$& $10.47$ & $5.9095e-11$& $10.47$ &  $6.4605e-11$& $10.47$ & $6.4568e-11$& $10.47$\\   
		&    $ 9^3\times 6 $		& $1.7452e-11$& $10.36$ & $1.7448e-11$& $10.36$ &  $1.9074e-11$& $10.35$ & $1.9060e-11$& $10.35$\\ 
        \addlinespace
        \multirow[t]{4}{*}{$ 10 $}&    $ 6^3\times 6 $ 		& $1.1866e-11$ &  & $1.8287e-11$&  & $1.2837e-11$&  & $1.5896e-11$& \\  
        &    $ 7^3\times 6 $ 		& $1.6818e-12$& $12.67$ & $2.5810e-12$& $12.70$ & $1.8153e-12$& $12.69$ & $2.2556e-12$& $12.67$\\  
        &    $ 8^3\times 6 $		& $3.1724e-13$& $12.49$ & $4.8556e-13$& $12.51$ & $3.4190e-13$& $12.50$ & $4.2581e-13$& $12.49$\\   
		&    $ 9^3\times 6 $		& $7.3856e-14$& $12.37$ & $1.1283e-13$& $12.39$ & $7.9528e-14$& $12.38$ & $9.9177e-14$& $12.37$\\ 
		\addlinespace 
		\multirow[t]{4}{*}{$ 11 $}&    $ 6^3\times 6 $ 		& $1.2702e-12$ &  & $6.3809e-12$&  & $1.2484e-13$&  &  $5.3666e-12$& \\  
        &    $ 7^3\times 6 $ 		& $1.8379e-13$& $12.54$ & $9.0907e-13$& $12.64$ & $1.2898e-14$& $14.73$ &  $7.6052e-13$& $12.68 $\\  
        &    $ 8^3\times 6 $		& $3.5127e-14$& $12.39$ & $1.7204e-13$& $12.47$ & $1.8509e-15$& $14.54$ & $1.4344e-13$& $12.49$\\ 
		&    $ 9^3\times 6 $		& $8.2625e-15$& $12.29$ & $4.0145e-14$& $12.36$ & $3.2092e-16$& $14.87$ & $3.3387e-14$& $12.38$\\    
		\addlinespace
		\multirow[t]{4}{*}{$ 12 $}&    $ 6^3\times 6 $ 		& $1.6497e-14$ &  & $1.0907e-13$&  & $1.2484e-13$&  & $1.3055e-14$&  \\  
        &    $ 7^3\times 6 $ 		& $1.6983e-15$& $14.75$ & $1.1286e-14$& $14.72$ & $1.2898e-14$& $14.73$ & $1.3540e-15$& $14.70$\\  
        &    $ 8^3\times 6 $		& $2.4633e-16$& $14.46$ & $1.6280e-15$& $14.50$ & $1.8509e-15$& $14.54$ & $1.9082e-16$& $14.67$ \\   
		&    $ 9^3\times 6 $		& $4.7705e-17$& $13.94$ & $3.0184e-16$& $14.30$& $3.2092e-16$& $14.88$ & $3.2960e-17$& $14.91$\\  
        \bottomrule
    \end{tabular}
\end{table*}\cref{tab:grid conv} shows the quadrature errors due to approximation of the integral in \cref{eq:integ smooth} and the convergence rates obtained with various rules of degrees 8 through 12. It shows that the convergence rates for the even and odd degree quadrature rules are close to $p+2$ and $p+1$, respectively. These rates are in agreement with those reported in \cite{chuluunbaatar2022pi,williams2020family} for quadrature rules on the tetrahedron and four-dimensional simplex. All the quadrature rules, except the degree 8 and 12 rules of \cite{chuluunbaatar2022pi}, the degree 11 line-LG rule, and the new degree 12 rule, achieve comparable accuracy levels and convergence rates. As noted in \cref{sec:efficiency line-LG}, the $q=11$ line-LG rule requires the same number of nodes as the $q=12$ line-LG rules; hence, identical rules are used for both cases. The new degree $12$ quadrature rule of this work and that of \cite{chuluunbaatar2022pi} seem to be significantly more accurate than the same degree quadrature rule in \cite{jaskowiec2021high} and the $q=12$ line-LG quadrature rule. They also achieve better accuracy for the highly oscillatory test case, as illustrated in \cref{fig:error tet}. On the other hand, despite the results in \cref{fig:error tet}, the $q=10$ rule in \cite{jaskowiec2021high} does not produce better results than the line-LG and the new quadrature rules when approximating \cref{eq:integ smooth}, suggesting that the improvement observed when approximating \cref{eq:integ J} is specific to the problem or mesh considered. It is not immediately clear why some of the same degree quadrature rules yield significantly different accuracy levels. Investigating ways to improve the accuracy of quadrature rules and operators that are derived using them, \eg, SBP derivative operators, can lead to drastic improvements in efficiency of numerical methods for PDEs, and is left for future studies.

\section{Conclusions}\label{sec:conclusions}
The main contribution of this paper is the development of fully-symmetric PI quadrature rules of very high orders. Using a novel approach for generating initial guesses and applying a node elimination strategy that takes into account the symmetry orbit structures of the optimal quadrature rule estimates, we extend the available set of symmetric PI quadrature rules from degree $50$ to $84$ on triangles, and from degree $20$ to $40$ on tetrahedra. Furthermore, the base quadrature rules before node elimination, the line-LG quadrature rules, are efficient in the sense that at very high degrees they require only about $10\%-20\%$ and $40\% -60\%$ more nodes than estimates of the optimal number of quadrature nodes on triangles and tetrahedra, respectively. Applying node elimination, we derived efficient quadrature rules, most of which have a number of nodes not more than $5\%$ and $25\%$ greater than the lower-bound estimates on triangles and tetrahedra, respectively. Numerical experiments show that the new quadrature rules yield accurate results for highly oscillatory functions and achieve the expected rate of convergence.

\appendix
\section{Supplementary Materials}
\label{app1}
The supplementary material associated with this paper can be found online at \url{https://github.com/OptimalDesignLab/SummationByParts.jl/tree/master/pi_quadrature_data}. It contains both the line-LG quadrature rules and the quadrature rules obtained after application of node elimination.




\bibliographystyle{elsarticle-num}
\bibliography{references}

\begin{thebibliography}{10}
\expandafter\ifx\csname url\endcsname\relax
  \def\url#1{\texttt{#1}}\fi
\expandafter\ifx\csname urlprefix\endcsname\relax\def\urlprefix{URL }\fi
\expandafter\ifx\csname href\endcsname\relax
  \def\href#1#2{#2} \def\path#1{#1}\fi

\bibitem{dunavant1985high}
D.~Dunavant, High degree efficient symmetrical {Gaussian} quadrature rules for
  the triangle, International Journal for Numerical Methods in Engineering
  21~(6) (1985) 1129--1148.

\bibitem{zhang2009set}
L.~Zhang, T.~Cui, H.~Liu, A set of symmetric quadrature rules on triangles and
  tetrahedra, Journal of Computational Mathematics (2009) 89--96.

\bibitem{xiao2010numerical}
H.~Xiao, Z.~Gimbutas, A numerical algorithm for the construction of efficient
  quadrature rules in two and higher dimensions, Computers \& Mathematics with
  Applications 59~(2) (2010) 663--676.

\bibitem{witherden2015identification}
F.~Witherden, P.~Vincent, On the identification of symmetric quadrature rules
  for finite element methods, Computers \& Mathematics with Applications
  69~(10) (2015) 1232--1241.

\bibitem{jaskowiec2021high}
J.~Ja{\'s}kowiec, N.~Sukumar, High-order symmetric cubature rules for
  tetrahedra and pyramids, International Journal for Numerical Methods in
  Engineering 122~(1) (2021) 148--171.

\bibitem{chuluunbaatar2022pi}
G.~Chuluunbaatar, O.~Chuluunbaatar, A.~Gusev, S.~Vinitsky, {PI}-type fully
  symmetric quadrature rules on the {$3-,\dots, 6-$}simplexes, Computers \&
  Mathematics with Applications 124 (2022) 89--97.

\bibitem{kirby2003aliasing}
R.~M. Kirby, G.~E. Karniadakis, De-aliasing on non-uniform grids: algorithms
  and applications, Journal of Computational Physics 191~(1) (2003) 249--264.

\bibitem{persson2009discontinuous}
P.-O. Persson, J.~Bonet, J.~Peraire, Discontinuous {Galerkin} solution of the
  {Navier--Stokes} equations on deformable domains, Computer Methods in Applied
  Mechanics and Engineering 198~(17-20) (2009) 1585--1595.

\bibitem{williams2019analysis}
D.~M. Williams, An analysis of discontinuous {Galerkin} methods for the
  compressible {Euler} equations: entropy and {$ L_2 $} stability, Numerische
  Mathematik 141~(4) (2019) 1079--1120.

\bibitem{fernandez2014review}
D.~C. Del Rey~Fern{\'a}ndez, J.~E. Hicken, D.~W. Zingg, Review of
  summation-by-parts operators with simultaneous approximation terms for the
  numerical solution of partial differential equations, Computers \& Fluids 95
  (2014) 171--196.

\bibitem{svard2014review}
M.~Sv{\"a}rd, J.~Nordstr{\"o}m, Review of summation-by-parts schemes for
  initial--boundary-value problems, Journal of Computational Physics 268 (2014)
  17--38.

\bibitem{hicken2016multidimensional}
J.~E. Hicken, D.~C. Del Rey~Fern{\'a}ndez, D.~W. Zingg, Multidimensional
  summation-by-parts operators: General theory and application to simplex
  elements, SIAM Journal on Scientific Computing 38~(4) (2016) A1935--A1958.

\bibitem{fernandez2018simultaneous}
D.~C. Del Rey~Fern{\'a}ndez, J.~E. Hicken, D.~W. Zingg, Simultaneous
  approximation terms for multi-dimensional summation-by-parts operators,
  Journal of Scientific Computing 75~(1) (2018) 83--110.

\bibitem{witherden2014analysis}
F.~D. Witherden, P.~E. Vincent, An analysis of solution point coordinates for
  flux reconstruction schemes on triangular elements, Journal of Scientific
  Computing 61 (2014) 398--423.

\bibitem{williams2020family}
D.~M. Williams, C.~V. Frontin, E.~A. Miller, D.~L. Darmofal, A family of
  symmetric, optimized quadrature rules for pentatopes, Computers \&
  Mathematics with Applications 80~(5) (2020) 1405--1420.

\bibitem{stroud1966gaussian}
A.~H. Stroud, D.~Secrest, Gaussian Quadrature Formulas, Prentice-Hall Englewood
  Cliffs, NJ, 1966.

\bibitem{worku2024quadrature}
Z.~A. Worku, J.~E. Hicken, D.~W. Zingg, Quadrature rules on triangles and
  tetrahedra for multidimensional summation-by-parts operators, Journal of
  Scientific Computing (Accepted), ArXiv Preprint arXiv:2311.15576 (2024).

\bibitem{solin2003higher}
P.~Solin, K.~Segeth, I.~Dolezel, Higher-order finite element methods, Chapman
  and Hall/CRC, 2003.

\bibitem{jaskowiec2020high}
J.~Ja{\'s}kowiec, N.~Sukumar, High-order cubature rules for tetrahedra,
  International Journal for Numerical Methods in Engineering 121~(11) (2020)
  2418--2436.

\bibitem{felippa2004compendium}
C.~A. Felippa, A compendium of {FEM} integration formulas for symbolic work,
  Engineering Computations 21~(8) (2004) 867--890.

\bibitem{proriol1957family}
J.~Proriol, Sur une famille de polynomes {\`a} deux variables orthogonaux dans
  un triangle, Comptes Rendus Hebdomadaires des S{\'e}ances de l'Acad{\'e}mie
  des Sciences 245~(26) (1957) 2459--2461.

\bibitem{koornwinder1975two}
T.~Koornwinder, Two-variable analogues of the classical orthogonal polynomials,
  in: R.~A. Askey (Ed.), Theory and Application of Special Functions, Academic
  Press, 1975, pp. 435--495.

\bibitem{dubiner1991spectral}
M.~Dubiner, Spectral methods on triangles and other domains, Journal of
  Scientific Computing 6 (1991) 345--390.

\bibitem{lyness1975moderate}
J.~N. Lyness, D.~Jespersen, Moderate degree symmetric quadrature rules for the
  triangle, IMA Journal of Applied Mathematics 15~(1) (1975) 19--32.

\bibitem{wang2023explicit}
W.~Wang, S.-A. Papanicolopulos, Explicit consistency conditions for fully
  symmetric cubature on the tetrahedron, Engineering with Computers 39~(6)
  (2023) 4013--4024.

\bibitem{chen2017entropy}
T.~Chen, C.-W. Shu, Entropy stable high order discontinuous {Galerkin} methods
  with suitable quadrature rules for hyperbolic conservation laws, Journal of
  Computational Physics 345 (2017) 427--461.

\bibitem{levenberg1944method}
K.~Levenberg, A method for the solution of certain non-linear problems in least
  squares, Quarterly of Applied Mathematics 2~(2) (1944) 164--168.

\bibitem{marquardt1963algorithm}
D.~W. Marquardt, An algorithm for least-squares estimation of nonlinear
  parameters, Journal of the Society for Industrial and Applied Mathematics
  11~(2) (1963) 431--441.

\bibitem{hicken2023summationbyparts}
J.~E. Hicken, Z.~A. Worku,
  \href{https://github.com/OptimalDesignLab/SummationByParts.jl}{{SummationByParts.jl,
  v0.2.1}}.
\newline\urlprefix\url{https://github.com/OptimalDesignLab/SummationByParts.jl}

\end{thebibliography}
\addcontentsline{toc}{section}{\refname}
\end{document}